\documentclass[twoside, 11pt]{article}
\usepackage{indentfirst,bm,latexsym,amsfonts,amssymb,amsmath,mathrsfs}
\RequirePackage{ifthen,calc}
\usepackage{theorem}
\usepackage{color}

\usepackage[colorlinks,citecolor=blue,linkcolor=blue]{hyperref}

 \def\beqlb{\begin{eqnarray}}\def\eeqlb{\end{eqnarray}}
 \def\beqnn{\begin{eqnarray*}}\def\eeqnn{\end{eqnarray*}}
\def\ra{\rightarrow}

 \numberwithin{equation}{section}
 \renewcommand{\theequation}{\arabic{section}.\arabic{equation}}
 \makeatletter
 \@addtoreset{equation}{section}
 \makeatother
 \newtheorem{theorem}{Theorem}[section]
 \newtheorem{definition}{Definition}[section]
 \newtheorem{hypothesis}{Hypothesis}[section]
 \newtheorem{lemma}{Lemma}[section]
 \newtheorem{proposition}{Proposition}[section]
 \newtheorem{corollary}{Corollary}[section]
 
 \newtheorem{example}{Example}[section]

 \def\qed{\hfill$\Box$\medskip}

 \def\bfE{{\mbox{\bf E}}}\def\bfP{{\mbox{\bf P}}}

\newcommand{\bcen}{\begin{center}}
\newcommand{\ecen}{\end{center}}
\newcommand{\bgeqn}{\begin{equation}}
\newcommand{\edeqn}{\end{equation}}

\def\ee{\varepsilon}
\def\re{\epsilon}

\def\rf{\rfloor}
\def\lf{\lfloor}
\def\rc{\rceil}
\def\lc{\lceil}

\def\L{{\cal L}}

\def\1{{\mathbf{1}}}

\def\ra{\rightarrow}

\def\qed{\hfill$\Box$\medskip}

\def\no{\nonumber}
\def\bfE{{\mathbb{E}}}
\def\mbfE{{\mathbf{E}}}
\def\bfP{{\mathbb{P}}}
\def\mbfP{{\mathbf{P}}}

\def\bfP{{\mathbb{P}}}
\def\bfR{{\mathbb{R}}}

\def\bfN{\mathbb{N}}

\def\iy{\infty}

\begin{document}
\title{\bf Asymptotic behavior for the quenched survival probability of a supercritical branching random
walk in random environment with a barrier
}
\author{You Lv\thanks{Email: lvyou@dhu.edu.cn },~~~~~Wenming Hong\thanks{Corresponding author. Email: wmhong@bnu.edu.cn}
\\
\\ \small School of Mathematics and Statistics, Donghua University,
\\ \small Shanghai 201620, P. R. China.
\\
\\ \small School of Mathematical Sciences $\&$ Laboratory of Mathematics and Complex Systems,
\\ \small Beijing Normal University, Beijing 100875, China
}
\date{}
\maketitle
\noindent\textbf{Abstract}: We introduce a random barrier to a supercritical branching random walk in an i.i.d. random environment $\{\L_n\}$ indexed by time $n,$ i.e.,
in each generation, only the individuals born below the barrier can survive and reproduce.
At generation $n$ ($n\in\bfN$), the barrier is set as $\chi_n+\ee n,$ where $\{\chi_n\}$ is a random walk determined by the random environment. Lv \& Hong (2024) showed that for almost every $\L:=\{\L_n\},$
the quenched survival probability (denoted by $\varrho_{\L}(\ee)$) of the particles system will be 0 (resp., positive) when $\ee\leq 0$ (resp., $\ee>0$). In the present paper,
we prove that $\sqrt{\ee}\log\varrho_{\L}(\ee)$ will converge in Probability/ almost surely/ in $L^p$ to an explicit negative constant (depending on the environment) as $\ee\downarrow 0$ under some integrability conditions respectively. This result extends the scope of the result of Gantert et al. (2011) to the random environment case.
\smallskip


\smallskip


\smallskip

\makeatletter 
\renewcommand\theequation{\thesection.\arabic{equation}}
\@addtoreset{equation}{section}
\makeatother 

%
%
%
%


 \def\bhypothesis{\begin{hypothesis}\sl{}\def\ehypothesis{\end{hypothesis}}}
 \def\bdefinition{\begin{definition}\sl{}\def\edefinition{\end{definition}}}
 \def\blemma{\begin{lemma}\sl{}\def\elemma{\end{lemma}}}
 \def\bproposition{\begin{proposition}\sl{}\def\eproposition{\end{proposition}}}
 \def\btheorem{\begin{theorem}\sl{}\def\etheorem{\end{theorem}}}
 \def\bcorollary{\begin{corollary}\sl{}\def\ecorollary{\end{corollary}}}
 \def\bremark{\begin{Remark}\rm{}\def\eremark{\end{Remark}}}
 \def\bexample{\begin{example}\rm{}\def\eexample{\end{example}}}

\allowdisplaybreaks
\section{Introduction}\label{Intro}
\subsection{Model} We consider a branching random walk on $\mathbb{R}$ in a time-inhomogeneous i.i.d. random environment
(BRWre), which is an extension of the time-homogeneous branching random walk (BRW). For a BRW, the reproduction law (including displacement and branching) of each generation is determined by a common point process, while for a BRWre, the reproduction law of each generation is sampled independently according to a common distribution on the collection of the point processes on $\bfR.$ The mathematical definition is as follows. 

Let $(\Pi,\mathcal{F}_{\Pi})$ be a measurable space and $\Pi\subseteq\tilde{\Pi}:=\{\mathfrak{m}:\mathfrak{m}~\text{is~a~point~process~on~}\bfR\}.$ The random environment $\L$ is defined as an i.i.d. sequence of random
 variables $\{\L_1$,~$\L_2$,~$\cdots,\L_n,\cdots\}$, where $\L_1$ takes values in $\Pi$.
 Let $\mu$ be the law of $\L$, then we call the product space $(\Pi^{\bfN}, \mathcal{F}_{\Pi}^{^{\bigotimes}\bfN}, \mu)$ the \emph{environment space}. For any realization $L:=\{L_1$,~$L_2$,~$\cdots,L_n,\cdots\}$ of $\L$, a time-inhomogeneous branching random walk driven by the environment $L$ is a process constructed as follows.

(1)~At time $0,$ an initial particle $\phi$ in generation $0$ is located at the origin.

(2)~At time $1,$ the initial particle $\phi$ dies and gives birth to $N(\phi)$ children who form the first generation. These children are located at $\zeta_i(\phi), 1\leq i\leq N(\phi),$ where the distribution of the random vector $X(\phi):=(N(\phi), \zeta_1(\phi),\zeta_2(\phi),\ldots)$ is $L_1.$ 

(3)~ Similarly, at generation $n+1,$ every particle $u$ alive at generation $n$ dies and gives birth to $N(u)$ children. If we denote $\zeta_i(u), 1\leq i\leq N(u)$ the displacement of the children with respect to their parent $u$, then the distribution of $X(u):=(N(u), \zeta_1(u),\zeta_2(u),\cdots)$ is $L_{n+1}.$ We should emphasize that conditionally on any given environment $L,$ all particles in this system reproduce independently. 

Conditionally on $\L,$ we write $(\Gamma,\mathcal{F}_{\Gamma}, \mbfP_{\L})$ for the probability space under which the
time-inhomogeneous branching random walk is defined. The probability $\mbfP_{\L}$ is usually called a {\it quenched law}.  We define the probability $\mathbf{P}:=\mu\bigotimes\mbfP_{\L}$ on the product space $(\Pi^{\bfN}\times\Gamma,\mathcal{F}_{\Pi}^{^{\bigotimes}\bfN}\bigotimes\mathcal{F}_{\Gamma})$ such that
\begin{eqnarray}\label{APP}\mathbf{P}(F\times G)=\int_{\L\in F}\mbfP_{\L}(G)~d\nu(\L), ~~F\in \mathcal{F}_{\Pi}^{^{\bigotimes}\bfN}, ~G\in\mathcal{F}_{\Gamma}.\no\end{eqnarray}
The marginal distribution of probability $\mathbf{P}$ on $\Gamma$ is usually called an {\it annealed law}. 
Throughout this paper, we consider the case $F=\Pi^{\bfN}.$ Hence without confusion we also denote $\mathbf{P}$ the annealed law and abbreviate $\mathbf{P}(\Pi^{\bfN}\times G)$ to~$\mathbf{P}(G).$ Moreover, we write $\mbfE_{\L}$ and~$\mathbf{E}$ for the corresponding expectation of $\mbfP_{\L}$ and~$\mathbf{P}$ respectively. This model can also be described by point process; see Mallein and Mi{\l}o\'{s} \cite{MM2016}.

We write $\mathbf{T}$ for the (random) genealogical tree of the process. For a given particle $u\in\mathbf{T}$ we write $V(u)\in\bfR$ for the position of $u$ and $|u|$ for the generation at which $u$ is alive. Then $(\mathbf{T}, V, \mbfP_{\L}, \mathbf{P})$ is called the {\it branching random walk in the time-inhomogeneous random environment $\L$} (BRWre). This model was first introduced in Biggins and Kyprianou \cite{BK2004}. If there exists a point process $\iota\in\Pi$ such that $\mathbf{P}(\L_1=\iota)=1$ (thus $\mathbf{P}(\L_i=\iota)=1, \forall i\in\bfN^+:=\{1,2,\cdots,n,\cdots\}$), then we usually call the environment a {\it degenerate environment} 
 and the BRWre degenerates to a BRW.
\subsection{A review on the barrier problem of BRW}
In the present paper we consider the barrier problem of BRWre. Let us first review a series of previous results on the barrier problem of BRW. 

In order to answer some questions about parallel simulations studied in Lubachevsky et al. \cite{LSW1989,LSW1990}, the barrier problem of BRW was first introduced in Biggins et al. \cite{BLSW1991}.

The so-called ``barrier" is actually a function $\varphi:\bfN\ra\bfR.$ Any particle $u$ and its descendants will be removed if $V(u)>\varphi(|u|),$ namely, a particle in this system can survive only if all its ancestors and itself were born below the barrier. The following notations express the barrier problem mathematically. A partial order $>$ on the tree $\mathbf{T}$ is defined as $u>v$ if $v$ is an ancestor of $u$. We write $u\geq v$ if $u>v$ or $u$ is $v$. For any $i\leq |u|,$ we denote $u_i$ the ancestor of $u$ in generation $i.$ We define an {\it infinite path} $u_{\infty}$ through $\mathbf{T}$ as a sequence of particles $u_{\infty}:=(u_i,i\in\bfN)$ such that $u_0=\phi~(\text{the initial particle})$, $\forall i\in\bfN, u_{i+1}>u_{i}.$  Let $\mathbf{T}_n:=\{u\in\mathbf{T}:|u|=n\}$ be the set of particles of generation $n$ and $\mathbf{T}_\infty$ the collection of all infinite paths through $\mathbf{T}.$ Denote $$\mathcal{S}_0:=\{\exists u_{\infty}=(u_0,u_1,u_2, \ldots u_n, \ldots)\in \mathbf{T_{\infty}}, \forall i\in\bfN,  V(u_i)\leq \varphi(i)\}$$ the event that the system survive when the barrier $\varphi$ was imposed on the BRW.

When a barrier is addressed to a BRW, a natural question is to consider the influence from the barrier on the survival/extinction of the system. Therefore, the barrier problem of BRW was usually considered under the assumptions that the underlying branching process is supercritical.  
Another basic assumption for the barrier problem of BRW is that the minimal displacement at time $n$ (denoted by $m_n$) grows at linear speed, i.e., \begin{eqnarray}\label{firstorder}\exists r\in\bfR, ~~\lim\limits_{n\rightarrow\infty}\frac{m_n}{n}=r, ~~{\rm a.s.}\end{eqnarray}
The sufficient conditions for the above convergence can be referred to Hammersley \cite{H1974}, Kingman \cite{K1975} and Biggins \cite{B1976}.


In the rest of this section and Section 2, we restate the related results in pervious papers under our setting. Denote $\varphi(i):=rn+\ee i^{\alpha}.$ Under the independence between the branching and displacement and some other mild integrable conditions, Biggins et al. \cite{BLSW1991} showed that
$\bfP(\mathcal{S}_0)>0$\footnote{We always write $\bfP$ for the law of a random model without random environment and $\bfE$ for the corresponding expectation.} when $\ee>0, \alpha=1$ and $\bfP(\mathcal{S}_0)=0$ when $\ee<0, \alpha=1.$


As a refined version of the above conclusion, Jaffuel \cite{BJ2012} showed that there exists a positive constant $\ee_*$ such that $\bfP(\mathcal{S}_0)>0$ when $\alpha=\frac{1}{3}, \ee>\ee^*$ and $\bfP(\mathcal{S}_0)=0$ when $\alpha=\frac{1}{3}, \ee<\ee^*$. 
Moreover, the independence of branching and displacement is not necessary for the method in \cite{BJ2012} and the following several papers. If the associated random walk is in the domain of
attraction of an $\alpha^*$-stable law, $\alpha^*\in(1,2)$, Liu and Zhang \cite{LZ2019} obtained an analogue of Jaffuel's result above. The main difference is the critical order of the barrier is $\frac{1}{\alpha^*+1}$ rather than $\frac{1}{3}.$

If the barrier causes extinction, scholars pay attention to the extinction rate and the total progeny. The extinction rate was obtained in Jaffuel \cite{BJ2012} for the case $\alpha=\frac{1}{3}, \ee\in(0,\ee^*)$ and in A\"{\i}d\'{e}kon \& Jaffuel \cite{AJ2011} for the case $\ee=0$ and $\alpha=1, \ee<0.$ For the case $\ee=0$, the tails of total progeny and the global leftmost position were studied in A\"{\i}d\'{e}kon \cite{A2010} and Addario-Berry \& Broutin \cite{AB2011}. A counterpart topic for branching brownian motion has also got a lot of attention in the past fifty years, see Berestycki et al. \cite{BBS2014} and the references therein.


Gantert et al. \cite{GHS2011} and B\'{e}rard \& Gou\'{e}r\'{e} \cite{BG2011} investigated the asymptotic behavior of $\bfP(\mathcal{S}_0)$ as $\ee\downarrow 0$. They both showed that \begin{eqnarray}\label{fe1}\lim\limits_{\ee\downarrow 0}\sqrt{\ee}\log\bfP(\mathcal{S}_0)\ra \rho,\end{eqnarray} where $\rho$ is a negative constant. For a general setting of the BRW, \cite{GHS2011} obtained the conclusion by a probabilistic approach while \cite{BG2011} focused on a special  cesr of BRW and gave additional precision on $\bfP(\mathcal{S}_0)$ by constructing a convolution equation. Berestycki et al. \cite{BBS2011} gave the analogue of \eqref{fe1} for branching Brownian motion.

In the present paper, we consider the corresponding asymptotic behavior in \eqref{fe1} in the context of BRWre. In the next subsection, we introduce some related results on BRWre.

\subsection{A review of the barrier problem of BRWre} 
Let us introduce the (log-)Laplace transforms of random point process $\L_n$
 \begin{eqnarray}\label{nota0}\kappa_n(\theta):=\log\mbfE_{\L}\left(\sum^{N(u)}_{i=1}e^{-\theta \zeta_i(u)}\right),~|u|=n-1,~ \theta\in[0,+\infty).\end{eqnarray}
 A basic assumption of the present paper is
 \begin{eqnarray}\label{au0}\exists\bar{\theta}>0, ~~\kappa(\bar{\theta})<+\infty,~~\exists \vartheta\in(0,\bar{\theta}),~~\kappa(\vartheta)=\vartheta\kappa'(\vartheta), ~~\kappa(0)\in(0, +\infty),\end{eqnarray}
 where $\kappa(\theta):=\mbfE(\kappa_1(\theta))$ (and thus $\kappa(\theta)=\mbfE(\kappa_n(\theta))$ for each $n$). This assumption was often set in other papers studying on BRWre. The model BRWre was first introduced in \cite{BK2004}.
 Recall that $m_n:=\min_{u\in \mathbf{T}_n}V(u)$ presents the minimal displacement in generation $n$. \cite{HL2014} proved that there is a finite constant $r^*$ such that
 \begin{eqnarray}\label{FOBRWre}\lim\limits_{n\rightarrow\infty}\frac{m_n}{n}=r^*,~~ \mathbf{P}-{\rm a.s.}\end{eqnarray} Conclusions on the central limit theorem of the BRWre can be found in \cite{GLW2014} and \cite{GL2016}. Large and moderate deviation principles for the counting measure have been obtained in \cite{HL2014} and \cite{WH2017} respectively.
 What inspires us most is the second order of the asymptotic behavior of $m_n$ considered in \cite{MM2016}. They showed that there exists a constant $c$ such that 
 \begin{align}\label{SOBRWre}\frac{m_n+\vartheta^{-1}K_n}{\log n}\ra c, ~~~n\ra\iy,~~~\text{in~~Probability}~~ \mathbf{P},\end{align}
where 
\begin{align}\label{Kn}K_n:=\sum_{i=1}^n\kappa_i(\vartheta),~~~~ K_0:=0.\end{align}
That is, the trajectory of
$\{m_n\}_{n\in\bfN}$ is around the random walk $\{\vartheta^{-1}K_n\}_{n\in\bfN}$ (noting that $K_n$ is the partial sum of i.i.d. random variables under $\mbfP$) with a logarithmic correction. This phenomenon exhibits the difference between the second order of the minimal displacement of the time-homogeneous case and that of the random environment case. More specifically, for a BRW satisfying some mild conditions, $m_n$ locates in the $\log n$-neighborhood of $rn$ (see \cite{HS2009} and \cite{AR2009}), where $r$ is the one in \eqref{firstorder}. Taking into account the above difference and the barrier ``$rn+\ee i^{\alpha}$" set for BRW in Section 1.1, we set the barrier for BRWre as  $$\varphi_{\L}(i):=-\vartheta^{-1}K_i+\ee i^{\alpha}$$ (but not as the form $``-\vartheta^{-1}i\mbfE K_1+\ee i^{\alpha}"$). We see $\varphi_{\L}(i)$ is a random barrier and the randomness comes from $K_i$ totally. For this barrier, we have obtained a series results in \cite{lv3} and \cite{lv4}.
Denote \begin{align*}\varrho_{\L}(\ee,\alpha):=\mbfP_{\L}(\exists u_{\infty}:=(u_1,u_2, \ldots u_n, \ldots)\in \mathbf{T_{\infty}}, \forall i\in\bfN,  V(u_i)\leq -\vartheta^{-1}K_i+\ee i^{\alpha})\end{align*}
the quenched survival probability after we add the barrier $\varphi_{_\L}$
and \begin{align*}Y_n:=\sharp\{|u|=n:~\forall i\leq n, ~V(u_i)\leq \varphi_{\L}(i)\},\end{align*}
the number of the surviving particles in generation $n$.

In \cite{lv3}, under some assumptions on $\kappa_1(\vartheta)$, the authors have shown that

{\rm (c1)}~$\varrho_{\L}(\ee,\alpha)=0,~{\rm \mathbf{P}-a.s.}$ when~$\alpha=\frac{1}{3}, \ee<\ee_c$ or $\alpha<\frac{1}{3}$, where $\ee_c$ is a positive explicit constant. (Hence we see $\varrho_{\L}(\ee,\alpha)=0,~{\rm \mathbf{P}-a.s.}$ as long as $\ee\leq 0$.)

{\rm (c2)}~The extinction rate can be characterized as $n^{-1/3}\log\mbfP_{\L}(Y_n>0)\ra t,~{\rm \mathbf{P}-a.s.}$ when $\alpha=\frac{1}{3}, 0<\ee<\ee_c$ or $\alpha<\frac{1}{3}, \ee\geq 0,$ where $t$ is a negative constants depending on $\ee$ and $\alpha$.

{\rm (c3)}~$\varrho_{\L}(\ee,\alpha)>0,~{\rm \mathbf{P}-a.s.}$ when~$\alpha=\frac{1}{3}, \ee>\ee_c$ or $\alpha>\frac{1}{3}, \ee>0$.

In \cite{lv4}, we further give the sufficient conditions for the extinction rate in $L^p$, that is,

{\rm (c4)} $n^{-1/3}\log\mbfP_{\L}(Y_n>0)\ra t,~{\rm in}~ L^p ~(p\geq 1).$ Of course, the sufficient conditions depend on $p$.

Combining (c1) with (c3), one sees that $\varrho_{\L}(\ee,1)=0, \mbfP-{\rm a.s.}$ if $\ee\leq 0$ and $\varrho_{\L}(\ee,1)>0, \mbfP-{\rm a.s.}$ if $\ee>0.$ In the present paper, we will give sufficient conditions to ensure the convergence of $\sqrt{\ee}\log\varrho_{\L}(\ee,1)$ in Probability, almost surely, or in $L^p$ as $\ee\downarrow 0$, which extends the result in Gantert et al. \cite{GHS2011} to the random environment case. Comparing with our previous work in \cite{lv3} and \cite{lv4}, to solve the problem in the present paper, we are facing more difficulties and need more preparations (see Section \ref{ewew11} for details). 

In addition, we refer to \cite{Oz2024} (and the references [17,18] therein) for a topic that a branching Brownian motion with a random trap (a mechanism partly similar to the barrier in the present paper).

\subsection{Organization}
The organization of this paper can be summarized as follows. The basic assumptions, main results and a related topic are given in Section 2. We prove the main results from Section 3 to Section 9 and refer to Section 2.4 for a detailed organization of the proof. The main task in the proof is how to prove the forthcoming \eqref{tpm}. We divide the proof of \eqref{tpm} into two steps and each step needs some preparations. 



\section{Main result}
\subsection{Assumptions}
Recall the notation in \eqref{nota0} and the assumption \eqref{au0}. Now we list the assumptions for our main results.

\noindent\emph {{$\bf(A1)$}. 
There exist $\lambda_0>3,~\lambda_1>2,~\lambda_{2}>2$ such that
\begin{eqnarray}\label{au1}%
\mathbf{E}\left(|\kappa_1(\vartheta)-\vartheta \kappa'_1(\vartheta)|^{\lambda_0}\right)+\mathbf{E}\left(\left[\frac{\mbfE_{\L}\left(\sum_{i=1}^{N(\phi)}|\zeta_i(\phi)+\kappa'_1(\vartheta)|^{\lambda_2}e^{-\vartheta \zeta_i(\phi)}\right)}{\mbfE_{\L}\left(\sum_{i=1}^{N(\phi)}e^{-\vartheta \zeta_i(\phi)}\right)}\right]^{\lambda_1}\right)<+\infty.
\end{eqnarray}
{$\bf(A1+)$} \eqref{au1} holds for some $\lambda_0>6, \lambda_1>3$ and $\lambda_2>2$.\\
{$\bf(A1-)$} \eqref{au1} holds for some $\min(\lambda_0,\lambda_1,\lambda_2)>2$. \footnote{Obviously,${\bf(A1+)}\Rightarrow{\bf(A1)}\Rightarrow{\bf(A1-)}.$ ${\bf(A1-)}$ will be mentioned in Section 3.}\\
{$\bf(A2)$}. There exist $\lambda_3>3,~\lambda_{4}>0$ such that 
\begin{eqnarray}\label{au3} \mbfE(|\kappa_1(\vartheta+\lambda_4)|^{\lambda_3})+\mbfE(|\kappa_1(\vartheta)|^{\lambda_3})+\mbfE([\log^+\mbfE_{\L}(N(\phi) ^{2})]^{\lambda_3})<+\infty,\end{eqnarray}
where we agree $\log^+\cdot:=\log\max\{\cdot,1\}$ and $\log^-\cdot:=|\log\min\{\cdot,1\}|.$ }

Throughout the paper, we denote \begin{eqnarray}\label{sigma*}\sigma^2:=\mathbf{E}\left(\Big(\kappa_1(\vartheta)-\vartheta\kappa'_1(\vartheta)\Big)^2\right),~~\sigma^2_*:=\vartheta^2\mathbf{E}(\kappa''_1(\vartheta)).\end{eqnarray}
\emph{{$\bf(A3)$}. The $\lambda_0, \lambda_2$ in {$\bf(A1)$} satisfy $\frac{\sigma^2}{\sigma^2_*}< \frac{\lambda_2-2}{\lambda_0-2}.$ Moreover, there exists $\lambda_{5}< -1, \lambda_6>2, \lambda_7>0$ such that
\begin{eqnarray}\label{au4}
\mathbf{E}\left(\left|\log^-\mbfE_{\L}\left(\1_{N(\phi)\leq |\lambda_5|}\sum_{i=1}^{N(\phi)}\1_{\{\vartheta \zeta_i(\phi)+ \kappa_1(\vartheta)\in [\lambda_5, \lambda_5^{-1}]\}} \right)\right|^{\lambda_6}\right)<+\infty,\end{eqnarray}
\begin{eqnarray}\label{au5}
\mathbf{E}\left(\left|\log^-\mbfE_{\L}\left(\1_{N(\phi)\leq |\lambda_5|}\sum_{i=1}^{N(\phi)}\1_{\{\vartheta \zeta_i(\phi)+ \kappa_1(\vartheta)\in [|\lambda_5|^{-1}, |\lambda_5|]\}} \right)\right|^{\lambda_6}\right)<+\infty,\end{eqnarray}
\begin{eqnarray}\label{au6}
\mathbf{E}\left(\left|\log^-\mbfE_{\L}\left(\sum_{i=1}^{N(\phi)}\1_{\{\vartheta \zeta_i(\phi)+ \kappa_1(\vartheta)\leq|\lambda_5|\}} \right)\right|^{\lambda_7}\right)<+\infty,\end{eqnarray}and we can find constant $\nu_0$ such that
\begin{eqnarray}\label{au7} \nu_0>2,~~ \nu_0\in [\lambda_6-1,\lambda_6),~~ \min\left\{\frac{\lambda_0}{2},\lambda_1,\frac{\lambda_3}{2}\right\}>\frac{\lambda_6}{\lambda_6-\nu_0}.\end{eqnarray}}

We refer to \cite[Example 2.4]{lv3} and \cite[Example 1.5 (2)]{lv4} for two examples satisfying all assumptions above.
These two examples are of different types: the former has a continuous, unbounded law of displacement which
is independent of the branching law, the latter has a discrete, bounded law of displacement
which may depend on the branching law. In \cite{lv3} and \cite{lv4}, we have verified these two examples satisfying all assumptions therein. Even though the assumptions in the present paper are more than those in \cite{lv3} and \cite{lv4}, we omit the extra verification since it can be done easily.

These assumptions are almost the same as the assumptions in \cite{lv4} (i.e., the sufficient conditions for existence of $L^p$-extinction rate, see (c4) in Section 1.3),
except the following two differences. One is that we set stronger integrability assumptions on $N(\phi)$ and the $\log$- quenched probabilities in \eqref{au4} and \eqref{au5}. \footnote{For example, we require $\mbfE([\log^+\mbfE_{\L}(N(\phi) ^{2})]^{\lambda_3})<+\infty$ in \eqref{au3} instead of $``\exists \re>0$, $\mbfE([\log^+\mbfE_{\L}(N(\phi) ^{1+\re})]^{\lambda_3})<+\infty$" in \cite{lv4}.} 
The other difference is that we add condition \eqref{au6}. But we observe that $\lambda_7\geq \lambda_6$ due to \eqref{au5}. Hence condition \eqref{au6} and $\lambda_7$ are only for a sharp characterization rather than an essential assumption.
We also refer to \cite[Section 2]{lv3} and \cite[Section 1.3]{lv4} for some detailed explanations on assumptions \eqref{au0} and ({\bf A1})-({\bf A3}).

\subsection{Main results}
Let us first recall the main result in \cite{lv2}. For two independent standard Brownian motions $B$ and $W$ (under a probability denoted by $\bfP$) starting at 0 ($B_0=W_0=0$), \cite{lv2} showed that there exists a real function on $\bfR$ such that
\begin{align}\label{bmt}\forall a>0,~\beta\in\bfR,~\hat\gamma(\beta):=\lim\limits_{t\rightarrow+\infty}\frac{-4a^2\log \bfP(\forall_{s\leq t} |B_s-\beta W_s|\leq a|W)}{t},~~{\rm a.s.,~and~}~~L^1,\end{align}
where $\hat\gamma(0)=\frac{\pi^2}{2}$ and $\hat\gamma(\cdot)$ is strictly increased at $\bfR^+.$

Denote $$\varrho_{\L}(\varepsilon):=\varrho_{\L}(\varepsilon,1),~~~ \lambda_8:=\min(\lambda_3, \lambda_7),~~~ \gamma_\sigma:=\sigma_*^2\hat\gamma\left(\frac{\sigma^2}{\sigma_*^2}\right),~~~ \gamma:=-\sqrt{\frac{\gamma_\sigma}{\vartheta}}.$$ Obviously, we have $\lambda_8>2.$ 
The following theorem, in which we give the sufficient conditions for the convergence of $\sqrt{\ee}\log\varrho_{\L}(\varepsilon)$ in Probability/ $L^p$ / almost surely, is the main result in the present paper.

\begin{theorem}\label{main}
We assume \eqref{au0} holds. Assumptions {$\bf(A1)$}-{$\bf(A3)$} are partly or totally needed in the following conclusions (1)-(6).


{\rm (1)}~If {$\bf(A1)$} holds, then for any $x>0,$ $\lim_{\ee\downarrow 0}\mbfP(\sqrt{\ee}\log\varrho_{\L}(\varepsilon)>\gamma+x)=0.$

{\rm (2)}~If {$\bf(A1+)$} holds, then $$\varlimsup_{\ee\downarrow 0}\sqrt{\ee}\log\varrho_{\L}(\varepsilon)\leq \gamma,~~ {\rm \mbfP-a.s.}$$ and hence $\varlimsup_{\ee\downarrow 0}\mbfE(\sqrt{\ee}\log\varrho_{\L}(\varepsilon))\leq \gamma$ (by Fatou's lemma).


{\rm (3)}~If {$\bf(A1)$}, {$\bf(A2)$} and {$\bf(A3)$} hold, then for any $x>0,$ $$\lim_{\ee\downarrow 0}\mbfP(\sqrt{\ee}\log\varrho_{\L}(\varepsilon)<\gamma-x)=0.$$ More precisely, there exists $\nu>\nu_0$ such that for any given $\lambda\in(2,\lambda_8)$,\footnote{Note that \eqref{tpm} can be stated equivalently as $\varlimsup_{\ee\downarrow 0}\ee^{\frac{(2-\min(\lambda_8,1+\nu_0))(\min(\lambda,\nu_0)-2)}{2\min(\lambda,\nu_0)}}\mbfP(\sqrt{\ee}\log\varrho_{\L}(\varepsilon)<\gamma-x)<+\infty$ due to the choice of $\nu_0$.}
\begin{eqnarray}\label{tpm} \varlimsup_{\ee\downarrow 0}\ee^{\frac{(2-\min(\lambda_8,1+\nu))(\min(\lambda,\nu)-2)}{2\min(\lambda,\nu)}}\mbfP(\sqrt{\ee}\log\varrho_{\L}(\varepsilon)<\gamma-x)<+\infty.\end{eqnarray}

{\rm (4)}~If ${\bf(A1)}, {\bf(A2)}$ and ${\bf(A3)}$ hold with $(\min(\lambda_8,1+\nu_0)-2)(\min(\lambda_8,\nu_0)-2)>3\min(\lambda_8,\nu_0),$
then $\varliminf_{\ee\downarrow 0}\sqrt{\ee}\log\varrho_{\L}(\varepsilon)\geq \gamma, {\rm \mbfP-a.s.}$

{\rm (5)}~If ${\bf(A1)}, {\bf(A2)}$ and ${\bf(A3)}$ hold, then for any $p\in[0,\min(\lambda_8,\nu_0)-2),$ we have $$\varlimsup_{\ee\downarrow 0}\mbfE(|\sqrt{\ee}\log\varrho_{\L}(\varepsilon)|^p)<+\infty.$$

{\rm (6)}~If ${\bf(A1)}, {\bf(A2)}$ and ${\bf(A3)}$ hold, 
then for any $p\in[0,\min(\lambda_8,\nu_0)-2)$,
we have
\begin{eqnarray}\sqrt{\ee}\log\varrho_{\L}(\varepsilon)\ra \gamma,~~ \ee\downarrow 0,~~{\rm in}~~ L^p.\end{eqnarray}
\end{theorem}

Obvious, Theorem \ref{main} (1) and (3) together imply that $\sqrt{\ee}\log\varrho_{\L}(\varepsilon)$ converges to $\gamma$  in Probability $\mbfP$ when ${\bf(A1)}, {\bf(A2)}$ and ${\bf(A3)}$ hold;  Theorem \ref{main} (2) and (4) together mean the convergence holds in the sense of $\rm \mbfP-a.s.$ when ${\bf(A1)}, {\bf(A2)}$ and ${\bf(A3)}$ hold with $(\min(\lambda_8,1+\nu_0)-2)(\min(\lambda_8,\nu_0)-2)>3\min(\lambda_8,\nu_0).$

According to the statements in \cite[Theorem 1.2]{GHS2011}, with our notation, the explicit expression of $\rho$ in \eqref{fe1} should be written as $\sqrt{\frac{\pi^2\sigma^2_*}{2\vartheta}}.$ Note that a degenerate environment satisfying \eqref{au0} implies $\sigma^2=0$ and recall that $\hat\gamma(0)=\frac{\pi^2}{2}$ mentioned in the beginning of this subsection. Then we see Theorem \ref{main} is consistent with \cite[Theorem 1.2]{GHS2011} when the random environment is degenerate.

\cite[Remark 2.12]{lv3} told that under some special setting, a slightly weaker version of \eqref{au4} is a necessary condition to ensure that for any $\ee>0$, $\varrho_{\L}(\ee)>0, \mbfP-{\rm a.s.}$ and hence a necessary condition for Theorem \ref{main}. Though we are not sure whether \eqref{au5} is a necessary condition, (note that \eqref{au5} seems not a counterpart assumption of \eqref{au4} since the right side of the barrier is set to be the deadly zone,) it is worth reminding that \cite[Remark 4.1]{lv4}  explains that  both \eqref{au4} and \eqref{au5} are almost the necessary conditions for the associated walk (introduced in Section 3) to satisfying the $L^p$-small deviation principle. Moreover, as far as we know, the small deviation estimate of the associated walk is an essential tool in the progress of the barrier problem for various branching random walks. At last, we mention that when the random environment is degenerate, \eqref{au4} and \eqref{au5} hold naturally under the assumption \eqref{au0}, see \cite[Proposition 1.4]{lv4}.

\subsection{Prospect}
 A future work which has a close connection to Theorem \ref{main} is to consider whether certain kinds of $N$-BRWre presents Brunet-Derrida behavior.  $N$-BRWre can be defined by extending the fixed reproduction law of a $N$-BRW to a random environment case. $N-$BRW is a branching selection particle system on the real line. In this model the total
size of the population at time $n$ is limited by $Z(n).$ At each generation $n$, every individual
dies while reproducing independently, making children around their current position according
to a fixed point processes, but only the $Z(n)$
leftmost children survive to form the $(n+1)$-th generation. The $N$-BRW has been studied in \cite{BG2011, M2017, M2018} under different settings. When $Z(n)\equiv N$ and the reproduction law is composed of a binary branching and a bounded walk step, which is the classical setting for a $N$-BRW, B\'{e}rard and Gou\'{e}r\'{e} \cite{BG2011} showed that both the maximal displacement $m^N_{1,n}$ and the minimal displacement $m^N_{2,n}$ of the $N$-BRW have the asymptotic behaviour $\lim_{n\ra\infty}\frac{m^N_{1,n}}{n}=\lim_{n\ra\infty}\frac{m^N_{2,n}}{n}:=r_N\in\bfR$ and $\lim_{N\ra\infty}(\log N)^2(r-r_N)=\rho$, where $r$ is the one in \eqref{firstorder} and $\rho$ is the one in \eqref{fe1}.  An intuitive understanding of this result from the view of \eqref{fe1} can be founded in \cite[Section 6.2]{S2015}. In this result, the $(\log N)^{-2}$ order magnitude of the difference $r_N-r$ was conventionally called as the Brunet-Derrida behavior of the $N$-BRW.  The Brunet-Derrida behavior stems from physicists' observation on the slowdown phenomenon in the wave propagation of some perturbed (indexed by $N$) F-KPP like
equations, for which the solution to the equation has a wave speed that is
slower than the standard speed (i.e., the speed w.r.t. the standard F-KPP equation) by a difference of order $(\log N)^{-2}$ when $N\ra +\infty$ (see \cite{BHP1994,BD1997,BD2001}). That is, both the asymptotic velocity $r_N$ and the aforementioned wave speed of the solution converge at the same slow rate $(\log N)^{-2}.$ In our context, based on Theorem \ref{main}, we could introduce a model $N-$BRWre with the setting that for a BRWre, among all particles at generation $n$, only the $Z(n)$
leftmost children survive to form the $(n+1)$-th generation. Then we consider, firstly, the relationship between the randomness of the corresponding quenched asymptotic velocity $r_N$ (or maybe it is more proper to write it as $r_{N,\L}$) and a proper setting for $Z(n)$ (whether $Z(n)$ depends on the random environment), and secondly, the existence of Brunet-Derrida behavior of the corresponding $r_{N,\L},$ i.e.,
whether the decay rate of $r_{N,\L}-r^*$ is $(\gamma+o(1))(\log N)^{-2}$ in some sense, where $r^*$ is the one in \eqref{FOBRWre} and from \eqref{SOBRWre} one can see $r^*=-\vartheta^{-1}\kappa(\vartheta).$

\subsection{\bf The outline of the proof}\label{ewew11}
Since we have developed the Mogul'ski\u{\i}'s estimate (a small deviation principle for random walk ,see \cite{Mog1974}) into the random environment case in \cite{lv2}, the proof of (1) and (2), which only involves the upper bound of $\varrho_{\L}(\varepsilon)$, can be adapted from \cite[Section3]{GHS2011}. Hence we only give a sketch on the proof of (1) and (2) after we introduce the associated walk of the BRWre in Section 3. In Section 9, we will see that (4) and (6) can be deduced quickly from (3) and (5) respectively. Moreover, the method used to prove (5) is similar to the one used in the proof of (3). Therefore, in the rest of the paper we will mainly concentrate on the proof of Theorem \ref{main} (3), i.e., the statement in \eqref{tpm}, which is the most difficult and technical part in the present paper.

In fact, thanks to the methods used in \cite{AJ2011,GHS2011,BJ2012}, when we have got the random environment version of Mogul'ski\u{\i}'s estimate in \cite{lv2}, we are very close to get, not only Theorem \ref{main} (1) and (2), but also almost all the results in \cite{lv3} (see the review in Section 1.3 (c1)-(c3)). However, it is not enough to get Theorem \ref{main} (3).
Theorem \ref{main} (3) can be viewed as an extension of the lower bound part in \eqref{fe1}.  The rigorous mathematically approach to obtain \eqref{fe1} can be found in  \cite{GHS2011}, \cite{BG2011} and  \cite{M2017}. As we have mentioned in Section 1, the method used in \cite{GHS2011} is purely probabilistic---combining a measure change, first-second moment argument and the small deviation principle. This method works for the BRW under a general setting (i.e., only some mild assumptions required) and gives the first order of the decay rate of $\bfP(\mathcal{S}_0)$ as $\ee\downarrow 0.$ Under some stronger assumptions (assuming that
the branching is binary and the random walk steps are bounded), \cite{BG2011} solved this problem in an analytical way---characterizing the survival probability as the solution of a non-linear convolution
equation and obtained a more precise estimate of the decay rate, that is, $|\log\bfP(\mathcal{S}_0)-\rho/\sqrt{\ee}|\leq O(\log \sqrt{\ee}).$ \cite{M2017} used almost the same (probabilistic) tools and gave an alternative proof for \eqref{fe1}. Compared to the proofs of the other results (in \cite{BLSW1991}, \cite{AJ2011} and \cite{BJ2012}) on the barrier problem for BRW, the extra difficulty in the proof of \eqref{fe1} is how to construct an auxiliary supercritical G.W. process whose survival probability is lower than $\bfP(\mathcal{S}_0).$
Actually, the difference between the proofs in \cite{M2017} and \cite{GHS2011} is mainly derived from the constructions of the auxiliary supercritical G.W. processes. Since the way in \cite{M2017} can be modified more easily to deal with the random environment case, in the proof of Theorem \ref{main} (3), we borrow the idea in \cite{M2017} to construct an auxiliary supercritical branching process in random environment (BPre). The auxiliary BPre will be given in Section 4. Meanwhile, the new difficulties come in train.
The way to estimate the survival probability of the auxiliary G.W. process from below
in \cite{M2017} strongly depends on the time-homogeneous property of the G.W. process and the BRW.
However, the estimate for the quenched survival probability of our auxiliary BPre appears more complexity. To complete the estimate, we need a further research
on the quenched probability that the BRWre with a barrier such like $-\vartheta^{-1}K_i+O(n^{-\frac{2}{3}})$ survives until time $n$ (the asymptotic behavior of this probability has been studied in \cite{lv3} preliminarily), see Section 6. A lower deviation of another auxiliary BPre is also needed, see Section 7. An essential tool used to prove the results in Section 6 and 7 is a corollary of a strong approximation theorem, which is given in Section 5. Based on all the preparations above, we finally complete the proof of \eqref{tpm} in Section 8.

Of course, the many-to-one lemma (a kind of change in measure which couples the BRWre we consider with an associated walk) is also a key tool in our proof. All the procedures in Sections 4-8 could be carried out only if we introduce the associated walk by a time-inhomogeneous version of many-to-one lemma, which is the main task in Section 3. Actually, compared with the study of BRW, a huge challenge during the study of BRWre
is directly reflected by the associated walk. By the many-to-one lemma, 
the associated walk w.r.t. BRWre is a random walk with random environment in time (RWre, see the definition in the next section) instead of the time homogeneous random walk in the context of BRW. We have obtained some limit behaviors of RWre in \cite{lv2}-\cite{lv5}. In Section 3, we give a slightly refined version of \cite[Proposition 4.3]{lv4} and \cite[Remark 2.1]{lv5}, which will support the proofs in the sequel.

All in all, we first owe the method developed in the study of the BRW with a barrier (summarized in the second paragraph in this subsection) and the asymptotic behavior of the minimal displacement of a BRWre (summarized in Section 1.3), without these predecessors' contributions we could do nothing on the barrier problem of a BRWre. However, staring from the previous work, it is still a long way to get Theorem \ref{main} in the present paper. 
In addition to some new ideas specially designed to complete the proof of Theorem \ref{main}, many of our previous efforts in \cite{lv1}-\cite{lv5} are indispensable for the present paper (though the motivations of these papers are not the issue addressed in the present paper). For example, we will see that at least one result or technique in each paper from \cite{lv1} to \cite{lv5} contributes the proof of Proposition \ref{pnb1}, directly or indirectly.



\section{Preparation 1: Many-to-one Formula and some properties of the associated walk}
The many-to-one formula (a kind of change of probabilities that transferring the information of the paths in the random genealogical tree to a random walk) is an essential tool in the study of the branching random walks.  It can be traced down to the early works of Peyri\`{e}re \cite{P1974}~and Kahane and Peyri\`{e}re \cite{KP1976}.
A version of time-inhomogeneous many-to-one formula has been introduced in Mallein \cite{M2015a}. The time-inhomogeneous many-to-one formula can also be applied to the study of BRWre, see \cite{MM2016}. On the other hand, for the time-homogeneous case, a bivariate version of many-to-one formula can be found in \cite{GHS2011}.  To prove \eqref{tpm},  we need a time-inhomogeneous bivariate version of many-to-one formula, which has been introduced in \cite{lv3}. For the convenience of reading, we restate the version in \cite{lv3} here. 
Let $\tau_{n,\L}$ be a random probability measure on $\bfR\times\bfN$ such that for any $x\in\bfR, A\in\bfN,$ we have
\begin{eqnarray}\label{m1}\tau_{n,\L}\left((-\infty,x]\times[0,A]\right)=\frac{\mbfE_{\L}\left(1_{\{N(u)\leq A\}}\sum^{N(u)}_{i=1}1_{\{\zeta_i(u)\leq x\}}e^{-\vartheta \zeta_i(u)}\right)}{\mbfE_{\L}\left(\sum^{N(u)}_{i=1}e^{-\vartheta \zeta_i(u)}\right)},~~|u|=n-1,\end{eqnarray}
where $\vartheta$ has been introduced in \eqref{au0}.
Hence we can see that the randomness of $\tau_{n,\L}$ comes entirely from $\L_n.$ 
Under the quenched law $\mbfP_{\L},$ we introduce a series of independent two-dimensional random vectors $\{X_n, \xi_n\}_{n\in\bfN^+}$ whose distributions are $\{\tau_{n,\L}\}_{n\in\bfN^+}.$ Define \begin{eqnarray}\label{shift0}S_0:=0, ~S_n:=\sum_{i=1}^{n}X_i, ~\forall n\in\bfN^+.\end{eqnarray}
The following lemma gives the relationship between $\{(S_n, \xi_n), n\in\bfN^+\}$ and the BRWre.
\begin{lemma}\label{mto}
{\bf (Many-to-one, \cite[Lemma 3.1]{lv3})} Let $v$ be a particle in generation $k\in\bfN.$ For any~$n\in\bfN^+$, any positive sequence~$\{A_i\}_{i\in \bfN^+}$ and any measurable function~$f:\bfR^n\rightarrow [0, +\infty),$ we have
\begin{eqnarray}\label{mto0}
&&\frac{\mbfE_{\L}\left[\sum_{|u|=n+k,u>v}\1_{\{N(u_{k+i})\leq A_i,1\leq i\leq n\}}e^{-\vartheta (V(u)-V(v))}f(V(u_{k+i})-V(v),1\leq i\leq n)\right]}{\mbfE_{\L}\left[\sum_{|u|=n+k, u>v}e^{-\vartheta (V(u)-V(v))}
\right]}\nonumber
\\&&~~~~~~~~~~~~~~~~=\mbfE_{\L}\left[f(S_{k+i}-S_k,1\leq i\leq n)\1_{\{\xi_{k+i}\leq A_i,1\leq i\leq n\}}\right],~~~{\rm\mathbf{P}-a.s.} \end{eqnarray}
\end{lemma}
Define \begin{eqnarray}\label{mtoT}
T_n:=\vartheta S_n+K_{n},
\end{eqnarray} where $\vartheta$, $K_n$ have been defined in \eqref{au0} and \eqref{Kn}. We can see $\{T_n\}$ is exactly the model ``a random walk with a random environment $\L$ in time (RWre)" considered in \cite{lv2} \footnote{According to the definition of RWre in \cite{lv2}, a more exact depiction of $\{T_n\}$ should be ``a random walk with a random environment $\{\tau_{n,\L}\}_{n\in\bfN^+}$ in time." But note that $\{\tau_{n,\L}\}_{n\in\bfN^+}$ is totally determined by $\L$ hence we can also say ``with random environment $\L$".}. We usually call $\{T_n\}$ as the {\it associated RWre with respect to the BRWre} or the \emph{associated walk}.  Now we are ready to prove Theorem \ref{main} (1) and (2).

\noindent{\bf A sketch of the proof of Theorem \ref{main} (1) and (2)}
By Lemma \ref{mto}, we see
\begin{eqnarray}\label{TandV}\mbfE_{\L}((T_1-\mbfE_{\L}T_1)^2)=\vartheta^2\kappa''_1(\vartheta),~~ \mbfE_{\L}T_1=\kappa_1(\vartheta)-\vartheta\kappa'_1(\vartheta),~~ \mbfE T_1=0.\end{eqnarray}
\cite[Section 4.1]{lv3} has shown that under the assumptions \eqref{au0} and \eqref{au1}, The RWre $\{T_n\}$ defined in \eqref{mtoT} satisfies all the assumptions in \cite{lv2}. That is, $\{T_n\}$ satisfies the convergence in probability part (resp. the almost surely part) of \cite[Corollary 2]{lv2} (the small deviation principle for RWre) when ${\bf(A1)}$ holds (resp. ${\bf(A1+)}$ holds). Then we can prove Theorem \ref{main} (1) and (2) by repeating the proof in \cite[Section 3]{GHS2011} step by step with the following modifications. The barrier $\ee i$ therein should be replaced by $\ee i-\vartheta^{-1}K_i$ and we use \cite[Corollary 2]{lv2} to replace all the Mogul'ski\u{\i} estimate used in \cite[Section 3]{GHS2011}. \qed

\cite[Corollary 2]{lv2} and its corollaries in \cite[Section 4.2]{lv3} play an important role in our study of the BRWre with a barrier, see \cite{lv3,lv4}. 
However, they are not enough to provide the proof of Theorem \ref{main} (3)-(6). 
Fortunately, \cite[Proposition 4.3]{lv4} and \cite[Remark 2.1]{lv5} have shown more properties of $\{T_n\}$, which are also needed in the proof of Theorem \ref{main} (3)-(6). Now we give more general versions of them as the following two lemmas. In order to keep the main line clear, we suggest readers admitting Lemma \ref{T-ubl} and \ref{T-ipl} and skipping over their proofs temporarily.

\begin{lemma}\label{T-ubl} The constants $a_1,a_2,b_1,b_2,a'_1,a'_2,z$ satisfy that $~b_1<a_1\leq a_2<b_2,~ b_1\leq a'_1< a'_2\leq b_2, z>0.$ Under assumptions \eqref{au0}, $({\bf A1-})$, $({\bf A2})$ and $({\bf A3})$, there exists $r<\frac{1}{2}$ and $\nu>\nu_0$ {\rm(}both $r$ and $\nu$ are independent of $a_1,a_2,b_1,b_2,a'_1,a'_2$ and $z${\rm)} such that
\begin{eqnarray}\label{T-ub}\varlimsup\limits_{n\rightarrow +\infty}\mbfE\left(~\left|\inf\limits_{x\in[a_1 \sqrt{n}, a_2 \sqrt{n}]}\log \mbfP_\L
\left(\begin{split}\forall_{i\leq \lfloor zn\rfloor}~ T_{i}\in[b_1\sqrt{n}, b_2\sqrt{n}],\\\xi_i\leq e^{n^r}, T_{\lfloor zn\rfloor}\in [a'_1 \sqrt{n}, a'_2 \sqrt{n}] \end{split}\Bigg|T_{0}=x\right)\right|^{\nu}~\right)<+\infty.~~\end{eqnarray}
\end{lemma}
{\bf Proof of Lemma \ref{T-ubl}}: 
This lemma is an extension of \cite[Proposition 4.3]{lv4}. We abbreviate \cite[Proposition 4.3]{lv4} as $\star$ in the rest of the proof. Here we only give a sketch since the extension is not so profound and the proof of $\star$ is too long. Let us first recall $\star$. $\star$ tells that there exists $\nu>\nu_0$ such that for any $b>a>0,$ if $b_2=-b_1=b$, $a'_2=a_2=-a'_1=-a_1=a$, then $\varlimsup\limits_{n\rightarrow +\infty}\mbfE\left(\left|p_{\L}(T;n,1,a,b,a,\frac{1}{2})\right|^{\nu}\right)<+\infty,$ where
\begin{eqnarray}\label{pt1}p_{\L}(T;n,z,a,b,c,r):=\inf_{|x|\leq a\sqrt{n}}\mbfP_\L
\left(\forall_{i\leq \lf zn\rf, i\in\bfN}~|T_{i}|\leq b\sqrt{n},~|T_{\lfloor zn\rfloor}|\leq c\sqrt{n},~\xi_{i}\leq e^{n^r}\big|T_{0}=x\right).~~~~~~~\end{eqnarray}

For a better understanding of the coming sketch, we review the approach of $\star$ as follows. Let us abbreviate $p_{\L}(T;n,z,a,b,a,r)$ as $\mathcal{X}_n$ until the end of this paragraph. For any $n$, we construct an event $Q_n$ and choose a constant $d_n$. Note that
\begin{eqnarray}\label{pt1+}&&\mbfE\left(|\mathcal{X}_n|^{\nu}\right)\leq \sum_{m=0}^{+\infty}\mbfP\left(|\mathcal{X}_n|^{\nu}\geq m\right)\no
\\&\leq& \sum_{m=d_n}^{+\infty}\mbfP\left(|\mathcal{X}_n|^{\nu}\geq m\right)+d_n\mbfP\left(Q^c_n\right)+\sum_{m=0}^{d_n-1}\mbfP\left(|\mathcal{X}_n|^{\nu}\geq m,Q_n\right)\no
\\&:=&D_{1,n}+D_{2,n}+D_{3,n}.\end{eqnarray}
In \cite[Section 3.2]{lv4}, we divide the whole proof of $\star$ into three parts. The goal of the $i$-th part ($i=1,2,3$) is to prove $\varlimsup\limits_{n\rightarrow +\infty}D_{i,n}<+\infty$.

The sketch is carried out in four steps. First we only consider the case that $a_1, a_2, b_1, b_2, a'_1, a'_2$ in \eqref{T-ub} satisfy the relationship that $b_2=-b_1:=b$, $a_2=-a_1:=a, a'_2=-a'_1:=a'.$ In step $i$ ($i=1,2,3$), by adapting the $i$-th part of the proof in \cite[Section 3.2]{lv4}, we check the upper limit of the corresponding $D_{i,n}$ is still finite in the context of $r<\frac{1}{2}$ and $a'\neq a$.
At last, we throw away the restriction $b_2+b_1=a_2+a_1=a'_2+a'_1=0$ in step 4.

%

\emph{Step 1} Note that the assumptions in Lemma \ref{T-ubl} is almost the same as the assumptions in $\star$ except only one slight difference. The difference is replacing the assumption
$$\exists\lambda_5<-1,~~ \mathbf{E}\left(\left[\log^-\mbfE_{\L}\left(\1_{N(\phi)\leq |\lambda_5|}\sum_{i=1}^{N(\phi)}\1_{\{\vartheta \zeta_i(\phi)+ \kappa_1(\vartheta)\in [0, |\lambda_5|]\}} \right)\right]^{\lambda_6}\right)<+\infty$$ in $\star$  by a stronger assumption \eqref{au5}. This only replacement will help us to achieve the goal in this step.

Note that for $n$ large enough,
\begin{eqnarray}&&p_{\L}(T;n,1,a,b,a',r)\no
\\&\geq& \prod_{i=1}^{n}\inf_{|x|\leq \left(\frac{n-i+1}{n}a+\frac{i-1}{n}a'\right)\sqrt{n}}\mbfP_{\L}\left(|T_{i}|\leq\left(\frac{n-i}{n}a+\frac{i}{n}a'\right)\sqrt{n},\xi_{i}\leq |\lambda_5|\Big|T_{i-1}=x\right)~~~~~\no
\\&\geq&\prod_{i=1}^{\lfloor zn\rfloor}\min\left\{\mbfP_{\L}\left(T_{i}\in[|\lambda_5|^{-1},|\lambda_5|],\xi_{i}\leq |\lambda_5||T_{i-1}=0\right),\mbfP_{\L}\left(T_{i}\in[\lambda_5,\lambda_5^{-1}],\xi_{i}\leq |\lambda_5||T_{i-1}=0\right)\right\}\no
\\&\geq& \prod_{i=1}^{\lfloor zn\rfloor}\left(\mbfP_{\L}\left(T_{i}\in[|\lambda_5|^{-1},|\lambda_5|],\xi_{i}\leq |\lambda_5||T_{i-1}=0\right)\mbfP_{\L}\left(T_{i}\in[\lambda_5,\lambda_5^{-1}],\xi_{i}\leq |\lambda_5||T_{i-1}=0\right)\right).\no\end{eqnarray}
Then we obtain that for each $\nu\in(\nu_0,\lambda_6),$
\begin{eqnarray}\label{Tu}~\varlimsup_{n\ra \infty}\sum_{m\geq \lceil n^{\frac{\nu}{\lambda_6-\nu}}\rceil}
\mbfP\left(\left|\log p_{\L}(T;n,1,a,b,a',r)\right|^{\nu}\geq m\right)<+\infty\no\end{eqnarray} (see \cite[(3.15)]{lv4}) for the analogous result) by adapting the first part in \cite[Section 3.2]{lv4}.

\emph{Step 2} Recalling assumption \eqref{au7} we see that there exists $r<\frac{1}{2}, \nu>\nu_0$ such that $r\lambda_3>\frac{\lambda_6}{\lambda_6-\nu}$. 
Now we redefine $I_n$ in \cite[(3.21)]{lv4} by replacing $\sqrt{n}$ with $n^r,$ then the upper bound of $\mathbf{P}(I^c_n)$ in \cite[(3.26)]{lv4} will be $\mathbf{P}(I^c_n)<c_7n^{1-r\lambda_3}$ rather than $\mathbf{P}(I^c_n)<c_7n^{1-\frac{\lambda_3}{2}}.$ Since the $r$ we choose also satisfies that $r\lambda_3>\frac{\lambda_6}{\lambda_6-\nu}$, the change (from $\sqrt{n}$ to $n^r$) will not break the conclusion ``$\varlimsup_{n\ra\infty} d_n\mbfP\left(Q^c_n\right)<+\infty$" in part 2 of \cite[Section 3.2]{lv4}.

\emph{Step 3} With our notation, though \cite[(3.28)]{lv4} and \cite[(3.32)]{lv4} should be rewritten as $${p}_{\L}(T; n,1,a,b,a',r)\geq p_{\L}(T;n,1,a,b,a',+\infty)-\sum_{i=1}^n\mbfP_{\L}(\xi_i>e^{n^r})$$
and $$\label{Ilow}\sum_{i=1}^n\mbfP_{\L}(\xi_i>e^{n^r})\leq ne^{-\frac{\lambda^2_4n^{r}}{3(\lambda_4+\vartheta)}}, ~~\text{on}~~ I_n,$$
the quantity $n\exp\{-\frac{\lambda^2_4n^{r}}{3(\lambda_4+\vartheta)}\}$ is still of the smaller order than RHS of \cite[(3.46)]{lv4}, which is a lower bound of $p_{\L}(T;n,1,a,b,a,+\infty).$
Hence when $a'=a$, the change (from ${p}_{\L}(T; n,1,a,b,a,\frac{1}{2})$ to ${p}_{\L}(T; n,1,a,b,a,r)$) will not break the final conclusion (the display therein following \cite[(3.46)]{lv4}) in part 3 of \cite[Section 3.2]{lv4}. Then we only need to consider the case $a'\in(0,a)$ since the case $a'>a$ is trivial.
It is indeed a little more difficult to deal with the case $a'\in(0,a)$ (i.e., to find the lower bound of $p_{\L}(T;n,1,a,b,a',+\infty)$), but the trick to overcome the difficulty has been presented in \cite[Section 4.2]{lv4}. That is, we can verify that even though $a'\in(0,a)$, the order of $n\exp\{-\frac{\lambda^2_4n^{r}}{3(\lambda_4+\vartheta)}\}$ is smaller than that of $p_{\L}(T;n,1,a,b,a',+\infty)$ by combining the arguments in \cite[Section 4.2]{lv4} and the part 3 of \cite[Section 3.2]{lv4}. So far, we have explained that the proof in \cite[Section 3.2]{lv4} can be adapted to prove that
\begin{eqnarray}\label{dongyuan}\exists r<\frac{1}{2}~{\rm and}~\nu>\nu_0,~~{\rm s.t.}~~\varlimsup\limits_{n\rightarrow +\infty}\mbfE\left(\left|p_{\L}(T;n,z,a,b,a',r)\right|^{\nu}\right)<+\infty\end{eqnarray}
when $z=1$, and hence \eqref{dongyuan} holds for any $z>0$ due to the arbitrariness of $a,b$ and $a'.$

\emph{Step 4} At last, we consider the case that $a_1, a_2, b_1, b_2, a'_1, a'_2$ in \eqref{T-ub} do not satisfy the relationship $b_2+b_1=a_2+a_1=a'_2+a'_1=0.$ Without loss of generality, the space homogeneous property of $T$ ensure that we can assume that $a_2=-a_1:=a.$ Let $\bar{a}_{2,j}:=\frac{m-j}{m}a+\frac{j}{m}a'_2$ and $\bar{a}_{1,j}:=-\frac{m-j}{m}a+\frac{j}{m}a'_1.$ We choose $m$ large enough such that $$\forall j\leq m-1,~~ \Delta_j:=\max\{|\bar{a}_{2,j+1}-\bar{a}_{2,j}|,|\bar{a}_{1,j+1}-\bar{a}_{1,j}|\}<\frac{\bar{a}_{2,j}-\bar{a}_{1,j}}{2}:=\alpha_{j}.$$ 
Let $\hat{\Delta}_j:=\min\{|\bar{a}_{2,j}-b_2|,~|\bar{a}_{1,j}-b_1|\}.$ 
Then according to the Markov property of $T$, a lower bound of the probability in \eqref{T-ub} could be given of the form $\prod_{j=0}^{m-1}p_{\L}(T; n,z/m,\alpha_{j},\alpha_{j}+\hat{\Delta}_j,\alpha_{j}-\Delta_j,r).$ Note that \eqref{dongyuan} tells that $\varlimsup\limits_{n\rightarrow +\infty}\mbfE(~|\log p_{\L}(T; n,z/m,\alpha_{j},\alpha_{j}+\hat{\Delta}_j,\alpha_{j}-\Delta_j,r)|^\nu~)<+\infty$ holds for any $j\leq m,$ which completes the proof of Lemma \ref{T-ubl}.\qed

%
%
\begin{lemma}\label{T-ipl} The constants $a_1,a_2,b_1,b_2,a'_1,a'_2$ satisfy that $b_1<a_1\leq a_2<b_2,~ b_1\leq a'_1< a'_2\leq b_2.$ 
Constant sequences $\{a_{1,n,i}\},$ $\{a_{2,n,i}\},$ $\{b_{1,n,i}\},$ $\{b_{2,n,i}\},$ $\{a'_{1,n,i}\}$ and $\{a'_{2,n,i}\}$ satisfy that \begin{eqnarray}\label{T-ipau}\lim\limits_{n\ra+\infty}\frac{\sum_{j=1}^2\left(\max\limits_{i\leq \lfloor tn\rfloor}|a_{j,n,i}-a_j\sqrt{n}|+\max\limits_{i\leq \lfloor tn\rfloor}|b_{j,n,i}-b_j\sqrt{n}|+\max\limits_{i\leq \lfloor tn\rfloor}|a'_{j,n,i}-a'_j\sqrt{n}|\right)}{\sqrt{n}}=0.~~~\end{eqnarray} $B$ and $W$ are two independent standard Brownian motions. 
If assumptions \eqref{au0} and ${\bf (A1-)}$ hold, then
\begin{eqnarray}\label{T-ip}
\mathfrak{P}_{n,\L}&:=&\inf_{z\in\left[a_{1,n,0},a_{2,n,0}\right]}\mbfP_\L
\left(\begin{split}\forall i\leq \lfloor tn\rfloor,~ T_{i}\in\left[b_{1,n,i},b_{2,n,i}\right], T_{\lfloor tn\rfloor}\in \left[a'_{1,n,\lfloor tn\rfloor},a'_{2,n,\lfloor tn\rfloor}\right]\end{split}\Bigg|T_{0}=z\right)\no
\\&~&~\ra\inf_{z\in\left[\frac{a_1}{\sigma_*},\frac{a_2}{\sigma_*}\right]}\bfP
\left(\begin{split}\forall s\leq t, \sigma_*B_s+\sigma W_s\in\left[b_1,b_2\right], \\ \sigma_*B_t+\sigma W_t\in \left[a'_1,a'_2\right]\end{split}\Bigg|W,W_0=0,B_0=z\right)~~~~~
\end{eqnarray}
in distribution as $n\ra +\infty$.
\end{lemma}
{\bf Proof of Lemma \ref{T-ipl}} First, recalling Lemma \ref{mto} and the definitions of $\sigma$ and $\sigma_*$,  we see
 $$\sigma^2=\mbfE((\mbfE_\L T_1)^2),~~~\sigma_*^2=\mbfE((T_1-\mbfE_\L T_1)^2)$$
 by a derict calculation. We also note that by Lemma \ref{mto}, ${\bf (A1-)}$ implies that
 $$\exists \lambda_0>2,~\lambda_1>2,~\lambda_2>2,~\mbfE(|\mbfE_\L T_1|^{\lambda_0})+\mbfE((\mbfE_{\L}|T_1-\mbfE_\L T_1|^{\lambda_2})^{\lambda_1})<+\infty.$$
Then comparing \eqref{T-ip} with the conclusion in \cite[Remark 2.1]{lv5}, we see  
\cite[Remark 2.1]{lv5} implies the truth of \eqref{T-ip} when $a_{1,n,0}=a_{2,n,0}=a_1=a_2=0,$ $b_{1,n,i}=a'_{1,n,i}=b_1\sqrt{n}$, $b_{2,n,i}=a'_{2,n,i}=b_2\sqrt{n}$, $a'_1=b_1<0<b_2=a'_2$ and $\sigma_*=t=1$.\footnote{A tip to the comparison: one can choose the $a_n$ and $b_n$ in \cite[Theorem 2.1]{lv5} as $(\log n)^{-1}$ and $n\log n$ respectively.} \cite[Remark 2.1]{lv5} was proved by applying the theory of strong approximation. Hence the proof is not involved any expansion of the two-sides boundaries crossing probability for a Brownian motion by trigonometric series like \cite[(A.4)]{GHS2011}.  Therefore, compared with the conclusion in \cite[Remark 2.1]{lv5}, the generalization on the starting point (from $0$ to $[a_{1}\sqrt{n}, a_{2}\sqrt{n}]$) and termination (from $[b_1\sqrt{n}, b_2\sqrt{n}]$ to $[a'_{1}\sqrt{n}, a'_{2}\sqrt{n}]$) of the trajectory will not bring new difficulty.

The further generalization from $a_{j}\sqrt{n}, a_{j}\sqrt{n}, a'_{j}\sqrt{n}$ to $a_{j,n,i}, b_{j,n,i}, a'_{j,n,i}, j=1,2,$ is convincing through the following observation. Note that for any given $\ee>0$, \eqref{T-ipau} ensures that
$$\mathfrak{P}_{n,\L}\leq\inf_{z\in\left[(\frac{a_1}{\sqrt{t}}+\ee)\sqrt{\lfloor tn\rfloor},(\frac{a_2}{\sqrt{t}}-\ee)\sqrt{\lfloor tn\rfloor}\right]}\mbfP_\L
\left(\begin{split}\forall i\leq \lfloor tn\rfloor,~ \frac{T_i}{\sqrt{\lfloor tn\rfloor}}\in\left(\frac{b_1}{\sqrt{t}}-\ee,\frac{b_2}{\sqrt{t}}+\ee\right),\\ \frac{T_{\lfloor tn\rfloor}}{\sqrt{\lfloor tn\rfloor}}\in \left(\frac{a'_1}{\sqrt{t}}-\ee,\frac{a'_2}{\sqrt{t}}+\ee\right)\end{split}\Bigg|T_{0}=z\right)$$
and
\begin{eqnarray}&&\mathfrak{P}_{n,\L}\geq\inf_{z\in\left[(\frac{a_1}{\sqrt{t}}-\ee)\sqrt{\lfloor tn\rfloor},(\frac{a_2}{\sqrt{t}}+\ee)\sqrt{\lfloor tn\rfloor}\right]}\mbfP_\L
\left(\begin{split}\forall i\leq \lfloor tn\rfloor,~ \frac{T_i}{\sqrt{\lfloor tn\rfloor}}\in\left(\frac{b_1}{\sqrt{t}}+\ee,\frac{b_2}{\sqrt{t}}-\ee\right),\\ \frac{T_{\lfloor tn\rfloor}}{\sqrt{\lfloor tn\rfloor}}\in \left(\frac{a'_1}{\sqrt{t}}+\ee,\frac{a'_2}{\sqrt{t}}-\ee\right)\end{split}\Bigg|T_{0}=z\right)\no\end{eqnarray}
hold for $n$ large enough in the sense of {\rm $\mathbf{P}-$a.s.}
 From the analysis above, we can easily obtain Lemma \ref{T-ipl} by applying the method used in the proof of \cite[Theorem 2.1]{lv5}. 

\section{Proof of \eqref{tpm}: the first step}
As stated in Section 2.4, the main task in the proof of Theorem \ref{main} (3)-(6) is to show \eqref{tpm}.

We divide the proof into two steps. In step 1, borrowing the idea in Mallein \cite[Section 3.3]{M2018}, we find a proper lower bound for $\sqrt{\ee}\log\varrho_{\L}(\varepsilon)$ by constructing several auxiliary branching processes. In step 2 (see Section 8), \eqref{tpm} will be obtained by estimating the lower bound given in step 1. 

{\it step 1} For any $c>1,\ee>0$, we want to construct a branching process whose survival probability is less than $\varrho_{\L}(c\varepsilon).$

First of all, we prove that there exists a constant $a$ large enough such that \begin{eqnarray}\label{ma}m(a):=\mathbf{E}\left(\log\mbfE_{\L}\left(\sum_{|u|=1}\1_{\{V(u)\leq a-\vartheta^{-1}K_1\}}\right)\right)>0,\end{eqnarray}
which could be verified from the following two facts.

Fact 1: Note that $\sum_{|u|=1}\1_{\{V(u)\leq a-\vartheta^{-1}K_1\}}\uparrow \sum_{|u|=1}1$ as $a\ra +\infty$ and hence by monotone convergence theorem, we have
$$\mbfE_{\L}\left(\sum_{|u|=1}\1_{\{V(u)\leq a-\vartheta^{-1}K_1\}}\right)\ra \mbfE_{\L}\left(N(\phi)\right), ~~\mbfP-{\rm a.s.}$$

Fact 2: Note that for each $a>|\lambda_5|$,
\begin{eqnarray}&&\left|\log\mbfE_{\L}\left(\sum_{|u|=1}\1_{\{V(u)\leq a-\vartheta^{-1}K_1\}}\right)\right|\leq
\log^-\mbfE_{\L}\left(\sum_{|u|=1}\1_{\{V(u)+\vartheta^{-1}K_1\leq|\lambda_5|\}}\right)+\log^+\mbfE_{\L}\left(N(\phi)\right).\no
\end{eqnarray}
Recalling assumptions \eqref{au6} and $\kappa(0)=\mathbf{E}\left(\log\mbfE_{\L}(N(\phi))\right)\in(0,+\infty)$ in \eqref{au0}, we deduce that $\lim_{a\rightarrow+\infty} m(a)=\kappa(0)$ by dominate convergence theorem, which means the truth of \eqref{ma}.


Define a positive sequence $\{\ee_n\}$ satisfying \begin{eqnarray}\label{fhao}n=\left\lfloor(\varsigma+z_n)\varepsilon_n^{\frac{-3}{2}}\right\rfloor,~~ z_n:=\frac{(c-1)\varsigma \varepsilon_n}{a-c\varepsilon_n},~~ \varsigma\in\bfN\end{eqnarray} and the choice of $\varsigma$ will be given later.
In the rest of the proof, we also write $\ee_n, z_n$ as $\ee, z$ if no confusion may arise. Define $$a_{n,l}=(l-1)(\varsigma n+\left\lfloor zn\right\rfloor)+\varsigma n,~~ b_{n,l}=l(\varsigma n+\lfloor zn\rfloor),~~ l=1,2,3,\cdots.$$
For any given $\ee$, we construct a branching process $\{G^*_l(\ee), l\in\bfN\},$ such that
$G^*_0(\ee)\equiv 1, G^*_l(\ee)=\sum_{j=1}^{G^*_{l-1}(\ee)}\eta^*_j(l,\ee)$, where $\{\eta^*_j(l,\ee),j\in\bfN^+\}$ is a sequence of i.i.d. random variables taking values in $\bfN$ with the common distribution that for every $\tau\in \bfN$,
\begin{eqnarray*}
&&\mbfP_{\L}(\eta^*_1(l,\ee)=\tau)
\\&:=&\mbfP_\L\left(\sharp\left\{|u|=b_{n,l}, u>v,  \forall i\leq b_{n,l}-b_{n,l-1}, V^{v}(u_{|v|+i})\leq h(c\ee,i,b_{n,l-1})\right\}= \tau\Big||v|=b_{n,l-1}\right).\no
\end{eqnarray*}
where $h(x,i,m):=x i- \vartheta^{-1}(K_{m+i}-K_{m}), V^v(u):=V(u)-V(v)$.
From the above construction one can see that $\{G^*_l(\ee), l\in\bfN\}$ is a branching process with random environment and the quenched survival probability $\rho^*_\L(\ee)$ of $\{G^*_l(\ee), l\in\bfN\}$ satisfies that
$$\rho^*_\L(\ee)\leq \varrho_\L(\ee),~~~\mbfP-{\rm a.s.}$$
Now let us find a lower bound of $\rho^*_\L(\ee).$ Let
\begin{eqnarray}\label{mi}m_i:=\mbfE_{\L}\left(\sum_{j=1}^{N(u)}\1_{\{\zeta_j(u)\leq a-\vartheta^{-1}\kappa_{i}\}}\right),~~ |u|=i-1,\end{eqnarray}
and $$m(l,\varepsilon):=w^{-\lfloor z n\rfloor}\Pi_{i=a_{n,l}+1}^{b_{n,l}}m_i,~~\bar m(l,\varepsilon):= \max\{m(l,\varepsilon), 1\},$$
where $w\in(1,e^{m(a)}).$
From the definition of $\ee$ and $z$, we see $$c\ee i\geq \ee\varsigma n +(i-\varsigma n)a,~~~\forall i\in[\varsigma n,\varsigma n+\lf z n\rf].$$
Then by Markov property we have
\begin{eqnarray}\label{pr1}
&&\mbfP_\L\left(\sharp\left\{|u|=b_{n,l}, u>v,  \forall i\leq b_{n,l}-b_{n,l-1}, V^{v}(u_{|v|+i})\leq h(c\ee,i,b_{n,l-1})\right\}\geq \bar m(l,\varepsilon)\Big||v|=b_{n,l-1}\right)\no
\\&\geq& p_1(l,\varepsilon)\bar{p}_2(l,\varepsilon),
\end{eqnarray}
where $$p_1(l,\varepsilon):=\mbfP_\L\left(\exists |u|=a_{n,l}, u>v, \forall i\leq \varsigma n, V^v(u_{|v|+i})\leq h(\ee,i,b_{n,l-1})\Big||v|=b_{n,l-1}\right),$$
$$\overline{p}_2(l,\varepsilon):=\mbfP_\L\left(\sharp\left\{|u|=b_{n,l}, u>v, \forall i\leq \lf zn\rf, V^v(u_{|v|+i})\leq h(a,i,a_{n,l})\right\}\geq\bar m(l,\varepsilon)\Big||v|=a_{n,l}\right).~~~$$
For any $l,$ we define a (time-inhomogeneous) branching process $\{Z_i(l),i\in\bfN\}$ such that
\begin{eqnarray}\label{defzn}Z_0(l)\equiv1,~~~ Z_i(l)=\sum_{j=1}^{Z_{i-1}(l)}\eta_j(a_{n,l}+i),\end{eqnarray} where for any $i$, $\eta_1(i), \eta_2(i),...,\eta_j(i),...$ are i.i.d. random variables taking values in $\bfN$ such that
\begin{eqnarray}\mbfP_{\L}(\eta_1(i)=\tau)=\mbfP_{\L}\left(\sum_{m=1}^{N(u)}\1_{\{\zeta_m(u)\leq a-\vartheta^{-1}\kappa_{i}\}}=\tau\right),~~ |u|=i-1,~~ \forall \tau\in\bfN.\no\end{eqnarray}
When $|v|=a_{n,l},$ we observe that under $\mbfP_{\L}$,
\begin{eqnarray}
&&\sharp\left\{|u|=b_{n,l}, u>v, \forall i\leq \lf zn\rf, V^v(u_{|v|+i})\leq
h(a,i,a_{n,l})\right\}\no
\\&\geq&\sharp\left\{|u|=b_{n,l},u>v, \forall i\leq \lf zn\rf, V^v(u_{|v|+i})-V^v(u_{|v|+i-1})\leq a-\vartheta^{-1}\kappa_{|v|+i}(\vartheta)\right\}\no
\\&:=&\Psi_{n,l}.\no
\end{eqnarray}
Hence we have \begin{eqnarray}\label{p2}\overline{p}_2(l,\varepsilon)\geq \mbfP_{\L}(\Psi_{n,l}\geq \bar m(l,\varepsilon))= \mbfP_{\L}(Z_{\lf zn\rf}(l)\geq \bar m(l,\varepsilon)):=p_2(l,\varepsilon).\end{eqnarray}
Now for any given $\ee$, we need to construct another branching process $\{G_l(\ee), l\in\bfN\}$ such that $G_0(\ee)\equiv 1, G_l(\ee)=\sum_{j=1}^{G_{l-1}(\ee)}\eta_j(l,\ee)$, where $\{\eta_j(l,\ee),j\in\bfN^+\}$ is a sequence of independent random variables with the common distribution
\begin{eqnarray*}
&&\mbfP_{\L}(\eta_1(l,\ee)=\lc\bar m(l,\varepsilon)\rc)=p_1(l,\varepsilon)p_2(l,\varepsilon),~~~\mbfP_{\L}(\eta_1(l,\ee)=0)=1-p_1(l,\varepsilon)p_2(l,\varepsilon).
\end{eqnarray*}
By the definition of $\bar m(l,\varepsilon), p_1(l,\varepsilon)$ and $p_2(l,\varepsilon),$ we see $\{G_l(\varepsilon), l\in\bfN\}$ is a Branching process with an i.i.d. random environment. Moreover,
$$\mbfP_{\L}(\eta_1^*{(l,\ee)}\geq \tau)\geq\mbfP_{\L}(\eta_1{(l,\ee)}\geq \tau), ~~\forall \tau\in\bfN$$
due to the relationship $\bar p_2(l,\varepsilon)\geq p_2(l,\varepsilon)$ and \eqref{pr1}.
 Let $\rho_{\L}(\varepsilon)$ be the quenched survival probability of $\{G_l(\varepsilon), l\in\bfN\}.$
From the above analysis one sees that
$$\varrho_\L(c\ee)\geq\rho^*_\L(\ee)\geq \rho_\L(\ee),~~~\mbfP-{\rm a.s.}$$
Let $\{f_{l,\varepsilon}\}_{l\in\bfN}$ be the generating functions of $\{G_l(\varepsilon), l\in\bfN\},$ namely
\begin{eqnarray*}
f_{l,\varepsilon}(s)&=&\mbfE_{\L}\left(s^{G_{l}(\ee)}\big|G_{l-1}(\ee)=1\right)=1-p_1(l,\varepsilon)p_2(l,\varepsilon)+p_1(l,\varepsilon)p_2(l,\varepsilon)s^{\lc\bar m(l,\varepsilon)\rc},~~ s\in[0,1].
\end{eqnarray*}
From \cite[Theorem 1]{AA1975} we see that
\begin{eqnarray}\label{5.4.5}
\mbfP(G_l(\varepsilon)>0)&\geq& \left[\frac{1}{\Pi_{j=1}^{l}f'_{j,\varepsilon}(1)}+\sum_{i=1}^{l}\frac{f''_{i,\varepsilon}(1)}{ f'_{i,\varepsilon}(1)}\frac{1}{\Pi_{j=1}^{i}f'_{j,\varepsilon}(1)}\right]^{-1}.\no
\end{eqnarray}
Since $\frac{f''_{i,\varepsilon}(1)}{ f'_{i,\varepsilon}(1)}=\lc\bar m(i,\varepsilon)\rc-1$, we have
$\mbfP(G_l(\varepsilon)>0)\geq\left[\sum_{i=1}^{l}\frac{\lc\bar m(i,\varepsilon)\rc}{\Pi_{j=1}^{i}f'_{j,\varepsilon}(1)}\right]^{-1}$
and hence
\begin{eqnarray}\label{lowrho}\varrho_\L(c\ee)\geq \left[\sum_{i=1}^{+\infty}\frac{\lc\bar m(i,\varepsilon)\rc}{\Pi_{j=1}^{i}f'_{j,\varepsilon}(1)}\right]^{-1},~~~\mbfP-{\rm a.s.}\end{eqnarray}
So far we have found a lower bound for $\varrho_\L(c\ee)$, where $c$ can be taken any value strictly larger than $1$.~Therefore, 
to show $\varliminf_{\ee\downarrow 0}\sqrt{\ee}\log\varrho_\L(\ee)\geq \gamma$ in Probability or almost surely, we should focus on the upper bound (of the decay rate) of
$$\mbfP\left(\sqrt{\varepsilon}\log\left[\sum_{i=1}^{+\infty}\frac{\lc\bar m(i,\varepsilon)\rc}{\Pi_{j=1}^{i}f'_{j,\varepsilon}(1)}\right]>-\gamma+x\right)$$
for any given $x>0$ small enough. Now we see that it is necessary to study the limit behaviors of $p_1(l,\varepsilon), p_2(l,\varepsilon)$ and $\bar m(l,\varepsilon)$, which are the subjects in the next three sections.
\section{Preparation 2: A corollary of strong approximation}
A tool used frequently in the forthcoming proof is a corollary of the celebrated Sakhanenko's strong approximation theorem. 
Let~$V_1,V_2,\ldots,V_n,\ldots$ be a sequence of independent random variables satisfying $\forall j,~ \bfE(V_j)=0$ and $\bfE(V^2_j)<+\infty$. Denote $\mathcal{D}_k:=\sum_{i=1}^k\bfE(V^{2}_i).$ Introduce a random broken line $\mathcal{V}(s), s\in\bfR^+$ such that $\mathcal{V}(0)=0,$ $\mathcal{V}(\mathcal{D}_k)=\sum_{i=1}^{k}V_i, k\in\bfN^+$ and $\mathcal{V}(\cdot)$ is linear, continuous on each interval $[\mathcal{D}_{k-1},\mathcal{D}_{k}].$ 
The following two theorems are known as the Sakhanenko's strong approximation theorem with power moment and the Cs\"{o}rg\H{o}~and~R\'{e}v\'{e}sz's estimate for Brownian motion.

\noindent\emph{{\bf Theorem \uppercase\expandafter{\romannumeral1}~(\cite[Theorem 1]{Sak2006})}
For any $\beta\geq 2,$ there exists a standard Brownian motion $B$ such that
\begin{eqnarray}\label{sak06}\forall x>0,~~~ \bfP\left(\sup_{s\leq \mathcal{D}_n}\big|\mathcal{V}(s)-B_s\big|\geq 2C_0\beta x\right)\leq
\frac{\sum_{k=1}^n\bfE(|V_k|^{\beta})}{x^{\beta}},\end{eqnarray}}
where $C_{0}$ is an absolute constant. \footnote{We always write $\bfP$ for the law and $\bfE$ the corresponding
expectation when the issues do not involve the random environment, e.g. all the conclusions in this section and the forthcoming Lemma \ref{pnbl1}.}

\noindent\emph{{\bf Theorem \uppercase\expandafter{\romannumeral2}~(\cite[Lemma 1]{CR1979})}} For a standard Brownian motion $B$ and a constant $D_1>2$, there exists a constant $D_2\in(0,+\infty)$ (depending only on $D_1$) such that
$$\forall x>0,~~~t>0,~~~ \bfP\left(\sup_{0\leq s\leq t}|B_s|\geq x\right)\leq D_2e^{-\frac{x^2}{D_1t}}.$$
The following corollary, which can be deduced from the two theorems above, plays a key role in the proofs in Sections 6-9.
\begin{corollary}\label{strapp}
$\{V_{i,\ee}, i\in\bfN, \ee\in\bfR^+\}$ is a triangle array of random variables such that $\bfE(V_{i,\ee})=0$ and for any $\ee,$ $\{V_{i,\ee}, i\in\bfN\}$ is a sequence of i.i.d. random variables with $\varlimsup_{\ee\downarrow 0}\bfE(|V_{1,\ee}|^{\beta})<+\infty$ for some $\beta\geq 2$. $U_{\ee}$ is a random variable which is independent of  $\{V_{i,\ee}, i\in\bfN\}$ and $\varlimsup_{\ee\downarrow 0}\bfE(|U_{\ee}|^{\beta})<+\infty.$
For any $m>0, l\in\bfN^+,$ if we can find $\iota>2$ such that
\begin{eqnarray}\label{sta}\iota\left[\varlimsup_{\ee\downarrow 0}\bfE(|V_{1,\ee}|^{2})l+\varlimsup_{\ee\downarrow 0}\bfE(|U_{\ee}|^{2})\right]m^{-2}\log(lm^{-\beta})\geq-1,\end{eqnarray}
then for $\ee>0$ small enough, there exists a positive constant $C$ independent of $l, m$ and $\ee$ such that
\begin{eqnarray}\label{stares}\forall m>0, ~~l\in\bfN^+,~~~~\bfP\left(\max_{i\leq l}\left|\sum_{k=1}^iV_{k,\ee}+U_{\ee}\right|\geq m\right)\leq C\frac{l}{m^\beta}.\end{eqnarray}
\end{corollary}
{\bf Proof of Corollary \ref{strapp}} Denote $$v_{r,\ee}:=\bfE(|V_{1,\ee}|^{r}),~u_{r,\ee}:=\bfE(|U_{\ee}|^{r}),~v_{r}:=\varlimsup_{\ee\downarrow 0}\bfE(|V_{1,\ee}|^{r}),~u_{r}:=\varlimsup_{\ee\downarrow 0}\bfE(|U_{\ee}|^{r}).$$
Note that \eqref{sta} implies that for any given $m>0,~ l\in\bfN^+$, there exist $\iota_0\in(0,1)$ and $\iota_1>0$ (independent of $l$ and $m$) such that
\begin{eqnarray}\label{sta1}~\frac{2+\iota_0}{1-\iota_0}[l(v_2+\iota_1)+(u_2+\iota_1)]m^{-2}\log(lm^{-\beta})\geq-1.~\end{eqnarray}
Let $t:=\sqrt{1-\iota_0}$ and $D_1:=2+\iota_0.$ 
From Theorem \uppercase\expandafter{\romannumeral1} we can find a  Brownian motion $B$ such that
$$\mathcal{P}_1:=\bfP\left(\max_{i\leq l}\left|\left(\sum_{k=1}^iV_{k,\ee}+U_{\ee}\right)-B_{iv_{2,\ee}+u_{2,\ee}}\right|\geq (1-t)m\right)\leq \left(\frac{2C_0\beta}{1-t}\right)^{\beta}\frac{lv_{\beta,\ee}+u_{\beta,\ee}}{m^\beta}.$$
Theorem \uppercase\expandafter{\romannumeral2} tells that there exists $D_2$ depending on $D_1$ such that
$$\mathcal{P}_2:=\bfP\left(\max_{i\leq l}\left|B_{iv_{2,\ee}+u_{2,\ee}}\right|\geq tm\right)\leq D_2\exp\left\{-\frac{t^2m^2}{D_1(v_{2,\ee}l+u_{2,\ee})}\right\}.$$
Therefore,
$$\mathcal{P}_2\leq D_2\exp\left\{-\frac{t^2m^2}{D_1((v_{2}+\iota_1)l+u_{2}+\iota_1)}\right\}\leq D_2\frac{l}{m^\beta}~~{\rm and}~~\mathcal{P}_1\leq \left(\frac{2C_0\beta}{1-t}\right)^{\beta}\frac{2l(v_{\beta}+u_{\beta})}{m^\beta}$$
hold for $\ee>0$ small enough, where the first inequality follows from \eqref{sta1}. Denote $C:=2\left(\frac{2C_0\beta}{1-t}\right)^{\beta}(v_{\beta}+u_{\beta})+D_2,$ which is a constant independent of $l$, $m$ and $\ee$. Finally we see \eqref{stares} follows from $\bfP\left(\max_{i\leq l}\left|\sum_{k=1}^iV_{k,\ee}+U_{\ee}\right|\geq m\right)\leq\mathcal{P}_1+\mathcal{P}_2$. \qed

\section{Preparation 3: Estimate on $p(n,b)$}
Let $b>0$ and denote $$p(n,b):=\mbfP_\L\left(\exists |u|=n, \forall i\leq n, V(u_i)\leq \frac{bi}{n^{2/3}}-\frac{K_{i}}{\vartheta}\right).$$ In this section, we study some asymptotic behaviors of $p(n,b)$ as $n\ra +\infty.$ If $\L_1$ is degenerate, the estimate on the decay rate of the constants sequence $\{p(n,b)\}_{n\in\bfN}$ is an important step in the study of $N$-BRW and BRW with a barrier, see \cite{M2018,GHS2011}.
In this section, we also give the convergence rate of the normalized quenched small deviation probability of a RWre (see Lemma \ref{pnbl2}), which is a further explosion on the small deviation probability considered in \cite{lv2}.

Throughout this section, we suppose that all the assumptions in Theorem \ref{main} (3) hold.

\begin{proposition}\label{pnb1p}
For any $b>0$, we have $\varliminf_{n\ra\infty}\mbfE(n^{-1/3}\log p(n,b))\geq \frac{\gamma}{\sqrt{b}}.$
\end{proposition}
{\bf Proof of Proposition \ref{pnb1p}} We divide this proof into two parts.

\emph{Step 1} By the first-second moment argument used in \cite[lemma 4.6]{GHS2011} (or see the proof of \cite[Theorem 2.5 (1b)]{lv3} w.r.t. the random environment case), for any $A_n> 1$ and constant $d\in\left(\sqrt{\frac{\gamma_{\sigma}}{\vartheta^3 b}},+\infty\right)$, we have
\begin{eqnarray}\label{pnb1}p\left(n,b\right)\geq \frac{\mbfE_{\L}\left(e^{T_{n}}\1_{\left\{\forall 1\leq i\leq n, ~T_i\in\left[\frac{\vartheta bi}{n^{2/3}}-dn^{\frac{1}{3}},\frac{\vartheta bi}{n^{2/3}}\right],~\xi_i\leq A_{n}\right\}}\right)}{1+(A_n-1)\sum\limits_{j=1}^{n}e^{\vartheta\left(bn^{1/3}+dn^{1/3}-\frac{bj}{n^{2/3}}\right)}P_{j,n}},~~{\rm \mathbf{P}-a.s.},\end{eqnarray}
where $\{T_n\}, \{\xi_{n}\}$ are both introduced in Section 3
and
\begin{eqnarray}\label{pnbpjn}P_{j,n}:=\sup_{y\in\bfR}\mbfP_{\L}\left(\forall_{j\leq i\leq n},\frac{T_i}{n^{1/3}}\in\left[\frac{b\vartheta i}{n}-d\vartheta ,\frac{b\vartheta i}{n}\right]\Big|T_j=y\right).\end{eqnarray}
Note that $P_{n,n}\equiv 1$ and hence for any $b'<b$,
\begin{eqnarray}\label{pnb2}p\left(n,b\right)\geq \frac{e^{n^{1/3}\vartheta b'}P_n}{A_n\sum\limits_{j=0}^{n}e^{\vartheta\left(bn^{1/3}+dn^{1/3}-\frac{bj}{n^{2/3}}\right)}P_{j,n}},~~{\rm \mathbf{P}-a.s.},\end{eqnarray}
where \begin{eqnarray}\label{pnbpn}P_n:=\mbfP_{\L}\left(\forall 1\leq i\leq n, T_i\in\left[\frac{\vartheta bi}{n^{2/3}}-\vartheta dn^{\frac{1}{3}},\frac{\vartheta bi}{n^{2/3}}\right], T_n\geq \frac{\vartheta b'n}{n^{2/3}},~~\xi_i\leq A_{n}|T_0=0\right).\end{eqnarray}
 Thanks to \eqref{au7}, we can choose constant $r$ such that $r<\frac{1}{2}$ and $r\lambda_3>\frac{\lambda_6}{\lambda_6-\nu_0}$. In the rest of the proof, we take $$A_n:=\exp\{n^{\frac{2}{3}r}\}.$$
Since $p\left(n,b\right)$ is non-increasing about $n,$ it is enough to show the asymptotic behavior of a sub-sequence of $\{p\left(n,b\right),n\in\bfN\}.$ Now we only consider the case that $n=N,2N,...,kN,...$, where $N\in\bfN.$
~From the above analysis, for $n=kN$, it is true that
\begin{eqnarray}\label{pnb3}p\left(n,b\right)\geq \frac{P_n}{A_nke^{n^{1/3}\vartheta (b-b'+d)}\sum\limits_{l=0}^{N-1}e^{-\frac{\vartheta blk}{n^{2/3}}}P_{(l+1)k,Nk}},~~{\rm \mathbf{P}-a.s.}\end{eqnarray}
Note that ${\bf (A1)}$ and ${\bf (A3)}$ imply ${\bf (A1+)}$ due to the relationships in \eqref{au7}. Recall the assumption ${\bf (A1+)}$ and the relationship \eqref{TandV}. Applying (the convergence almost surely part of) \cite[Corollary 2]{lv2}, we have
$$\lim_{n\ra\infty} \frac{1}{n^{1/3}}\log P_{(l+1)k,Nk}= -\frac{N-l-1}{N}\frac{\gamma_{\sigma}}{\vartheta^2d^2}, ~{\rm \mathbf{P}-a.s.}$$
Recalling that $d>\sqrt{\frac{\gamma_{\sigma}}{\vartheta^3 b}}$, we see
\begin{eqnarray}\lim_{n\ra\infty} \frac{1}{n^{1/3}}\log\left(\sum\limits_{l=0}^{N-1}e^{-\frac{\vartheta blk}{n^{2/3}}}P_{(l+1)k,Nk}\right)&=& \max_{0\leq l\leq N-1}\left(-\frac{\vartheta bl}{N}-\frac{N-l-1}{N}\frac{\gamma_{\sigma}}{\vartheta^2d^2}\right)\no
\\&=& -\frac{N-1}{N}\frac{\gamma_{\sigma}}{\vartheta^2d^2}, ~~{\rm \mathbf{P}-a.s.}.\no\end{eqnarray}
Note that for any fixed $N$, $\sum_{l=0}^{N-1}e^{-\frac{\vartheta blk}{n^{2/3}}}P_{(l+1)k,Nk}\leq \sum_{l=0}^{N-1}e^{-\frac{\vartheta blk}{n^{2/3}}}\leq (1-e^{-\frac{\vartheta bn^{1/3}}{N}})^{-1}\leq 2$ as long as $n$ large enough. Then Fatou's lemma tells that 
\begin{eqnarray}\label{pnb4}\varlimsup_{n\ra\infty} \mbfE\left(n^{-1/3}\log\left(\sum\limits_{l=0}^{N-1}e^{-\frac{\vartheta blk}{n^{2/3}}}P_{(l+1)k,Nk}\right)\right)\leq -\frac{N-1}{N}\frac{\gamma_{\sigma}}{\vartheta^2d^2}.\end{eqnarray}
\emph{Step 2} Now we find a lower bound of the numerator in \eqref{pnb3}.
Choose $\delta\in\left(0,\min\{\frac{2}{3}\frac{\vartheta(b-b')}{|\lambda_5|},\frac{\vartheta d}{2|\lambda_5|}\}\right).$ Let
\begin{eqnarray}\label{last62}D\in\bfR^+,~~\delta_n:=\left\lf\delta n^{1/3}\right\rf,~~ D_n:=\lf Dn^{2/3}\rf,~~ D_{n,l}:=\delta_n+lD_n,~~ s_n:=\left\lf\frac{n-\delta_n}{D_n}\right\rf-1.\end{eqnarray}
Denote
\begin{eqnarray}\label{pnbnotat}I^*(t_1, t_2):=\left\{\forall t_1\leq i\leq t_2,~ T_i\in\left[\frac{\vartheta bi}{n^{2/3}}-\vartheta dn^{\frac{1}{3}}, \frac{\vartheta bi}{n^{2/3}} \right],\xi_i\leq A_{n}\right\},\end{eqnarray}
$$I(t_1, t_2):=\left\{\forall t_1\leq i\leq t_2,~ T_i\in\left[2D\vartheta b-\vartheta dn^{\frac{1}{3}}, 0\right],\xi_i\leq A_{n}\right\},$$
$$\left[c_1, c_2\right]_{n,l}:=\left[c_1+\frac{\vartheta blD_n}{n^{2/3}},c_2+\frac{\vartheta blD_n}{n^{2/3}}\right].$$
For $n$ large enough, Markov property tells that
\begin{eqnarray}\label{pnb>}
P_n&\geq&\prod_{m=1}^{\delta_n}\mbfP_\L
(T_{m}\in [\lambda_5,\lambda^{-1}_5],\xi_{m}\leq A_n|T_{m-1}=0)\nonumber
\\&\times&\prod_{l=0}^{s_n-1}\inf_{z\in\left[\lambda_5\delta_n, \lambda^{-1}_5\delta_n\right]_{n,l}}\mbfP_\L
\Big(I^*(D_{n,l},D_{n,l+1}), T_{D_{n,l+1}}\in \left[\lambda_5\delta_n, \lambda^{-1}_5\delta_n\right]_{n,l+1}\big|T_{D_{n,l}}=z\Big)\no
\\&\times&\inf_{z\in\left[\lambda_5\delta_n, \lambda^{-1}_5\delta_n\right]_{n,s_n}}\mbfP_\L
\Big(I^*(D_{n,s_n},n), T_{n}\geq \vartheta b'n^{1/3}\big|T_{D_{n,s_n}}=z\Big)\no
\\&:=&\prod_{m=1}^{\delta_n}q_m\prod_{l=0}^{s_n-1}\varphi^*_{n,l}\cdot\varphi^*_{n,end},
\end{eqnarray}
 where $\lambda_5< -1$ is the one in assumption \eqref{au4}.
According to the space-homogeneous property of the model and the fact that $\forall i\leq D_n, \frac{\vartheta bi}{n^{2/3}}\in[0,\vartheta bD],$ we see
$$\varphi^*_{n,l}\geq \varphi_{n,l}:=\inf_{z\in\left[\lambda_5\delta_n, \lambda^{-1}_5\delta_n\right]}\mbfP_\L
\left(I(D_{n,l},D_{n,l+1}), T_{D_{n,l+1}}\in\left[\vartheta bD+\lambda_5\delta_n, \lambda^{-1}_5\delta_n\right]\Big|T_{D_{n,l}}=z\right),$$
\begin{eqnarray}\label{last61}\varphi^*_{n,end}\geq \inf_{z\in\left[\lambda_5\delta_n, \lambda^{-1}_5\delta_n\right]}\mbfP_\L
\left(I(D_{n,s_n},n), T_{n}\geq \vartheta b'n^{1/3}-\frac{\vartheta bs_nD_n}{n^{2/3}}\Big|T_{D_{n,s_n}}=z\right).\end{eqnarray}
Since we have chosen that $\delta\in(0, \frac{2}{3}\frac{\vartheta(b'-b)}{|\lambda_5|}),$ we see
$\lambda_5\delta_n-(\vartheta b'n^{1/3}-\frac{\vartheta bs_nD_n}{n^{2/3}})\geq \frac{1}{4}\vartheta(b-b')n^{1/3}.$
Then for $n$ large enough, \eqref{last61} tells that
\begin{eqnarray}\varphi^*_{n,end}&\geq& \inf_{z\in\left[\lambda_5\delta_n, \lambda^{-1}_5\delta_n\right]}\mbfP_\L
\left(I(D_{n,s_n},n), T_{n}\geq \vartheta b'n^{1/3}-\frac{\vartheta bs_nD_n}{n^{2/3}}\Big|T_{D_{n,s_n}}=z\right)\no
\\&\geq&\mbfP_\L
\left(\forall D_{n,s_n}\leq i\leq D_{n,s_n}+\lfloor 2Dn^{2/3}\rfloor, |T_{i}|\leq \min\left\{\frac{\delta n^{1/3}}{2|\lambda_5|},\frac{\vartheta (b-b')n^{1/3}}{4}\right\}\Big|T_{D_{n,s_n}}=0\right)\no
\\&:=&\varphi_{n,end}.\end{eqnarray}
Moreover, we observe that for any given $n$, $\{\varphi_{n,l}\}_{0\leq l\leq s_n-1}$ is an i.i.d. sequence and the common distribution is the same as 
\begin{eqnarray}
\varphi_{n}:=\inf_{z\in\left[\lambda_5\delta_n, \lambda^{-1}_5\delta_n\right]}\mbfP_\L\left(I(0,D_{n}), T_{D_{n}}\in\left[\vartheta bD+\lambda_5\delta_n, \lambda^{-1}_5\delta_n\right]\Big|T_{0}=z\right).
\end{eqnarray}
From the above analysis and \eqref{pnb>} we have
\begin{eqnarray}\label{pnb0}&&\mbfE\left(n^{-1/3}\log\mbfP_{\L}\left(\forall 1\leq i\leq n, T_i\in\left[\frac{\vartheta bi}{n^{2/3}}-\vartheta dn^{\frac{1}{3}},\frac{\vartheta bi}{n^{2/3}}\right], T_n\geq \frac{\vartheta b'n}{n^{2/3}},~~\xi_i\leq A_{n}\right)\right)\no
\\&\geq&n^{-1/3}[\delta_n\mbfE\left(\log q_1\right)+s_n\mbfE\left(\log \varphi_n\right)+\mbfE\left(\log\varphi_{n,end}\right)].\end{eqnarray}
The following displays (from \eqref{pnbsh1} to \eqref{pnbresult}) show the asymptotic behavior of $\mbfE\left(\log \varphi_n\right).$ Let $$\tilde{\varphi}_{n}:=\inf_{z\in\left[-\lambda\delta_n, \lambda_5\delta_n\right]}\mbfP_\L
\Bigg(\forall i\leq D_{n},~ T_{i}\in\left[\vartheta bD-\vartheta dn^{\frac{1}{3}},0\right], T_{D_n}\in \left[\vartheta bD+\lambda_5\delta_n, \lambda^{-1}_5\delta_n\right]\Big|T_{0}=z\Bigg),$$
then it is true that \begin{eqnarray}\label{pnbsh1}~~~~~\tilde\varphi_{n}-\sum_{i=1}^{D_n}\mbfP_{\L}(\xi_i>A_n)\leq \varphi_{n}\leq \tilde\varphi_{n}, ~~{\rm \mathbf{P}-a.s.}\end{eqnarray}
Lemma \ref{T-ipl} tells that
\begin{eqnarray}\label{pnbsh2}\varphi_{n}\ra\varphi(D):=\inf_{z\in\left[\frac{\lambda_5\delta}{\sigma_*},\frac{\delta}{\lambda_5\sigma_*}\right]}\mbfP_\L
\left(\begin{split}\forall_{s\leq D},~ \sigma_*B_s+\sigma W_s\in\left[-\vartheta d,0\right],\\ \sigma_*B_D+\sigma W_D\in [\lambda_5\delta,\delta/\lambda_5]\end{split}\Bigg|W,W_0=0,B_0=z\right)~~\end{eqnarray}
in distribution, where $B$ and $W$ are two independent standard Brownian motions. Our task in the next few paragraphs is to check $\varphi_{n}, \varphi(D), \tilde\varphi_{n}$ satisfying all the assumptions in Lemma \ref{pnbl1}.

From the choice of $\delta$ we see $\lambda_5\delta-(2D\vartheta b-\vartheta dn^{\frac{1}{3}})>\frac{1}{3}\vartheta dn^{\frac{1}{3}}$ for $n$ large enough.
~Therefore, by Lemma \ref{T-ubl} we see that there exists $\nu>\nu_0(>2)$ such that
\begin{eqnarray}\label{pnbsh3} \varlimsup_{n\ra\infty}\mbfE(|\log\varphi_{n}|^{\nu}+|\log\varphi_{n,end}|^{\nu})<+\infty~\end{eqnarray}
and hence $\varlimsup_{n\ra\infty}\mbfE(|\log\tilde{\varphi}_{n}|^{\nu})<+\infty.$

Now we estimate the tail of $\sum_{i=1}^{D_n}\mbfP_{\L}(\xi_i>A_n).$ From \cite[the second display after (4.15)- (4.17)]{lv3} we see $$\mbfP_{\L}(\xi_i>A_n)\leq \exp\{-\lambda_3v_1n^{\frac{2}{3}r}\}\mbfE_{\L}(N(u)^{1+\lambda_4})^{v_1}e^{(1-v_1)\kappa_i(\vartheta+\lambda_4)-\kappa_i(\vartheta)},$$
where $|u|=i-1, v_1:=\frac{\lambda_4}{\vartheta+\lambda_4}$ and $\lambda_4$ is the one in assumption ({\bf A2}). Markov inequality and ({\bf A2}) ensure that $$\exists c_3,~~~\mbfP\left(\mbfE_{\L}(N(u)^{1+\lambda_4})^{v_1}e^{(1-v_1)\kappa_i(\vartheta+\lambda_4)-\kappa_i(\vartheta)}\geq \exp\left\{\frac{1}{2}\lambda_3v_1n^{\frac{2}{3}r}\right\}\right)\leq c_3n^{-\frac{2}{3}r\lambda_3}.$$
Recall that we have chosen $r\lambda_3>\frac{\lambda_6}{\lambda_6-\nu_0}$ and hence $r\lambda_3>1$. Hence for $n$ large enough, there exists $c_4>0$ such that
\begin{eqnarray}\label{pnbsh4}\mbfP\left(\sum_{i=1}^{D_n}\mbfP_{\L}(\xi_i>A_n)>\exp\left\{-\frac{1}{3}\lambda_3v_1n^{\frac{2}{3}r}\right\}\right)<c_4n^{\frac{2}{3}-\frac{2}{3}r\lambda_3}\ra 0.\end{eqnarray}
Moreover, \cite[Theorem 3.1]{lv1} tells that
\begin{eqnarray}\label{pnbsh4+}\mbfP\left(\log\varphi(D)\right)<+\infty.\end{eqnarray}
According to \eqref{pnbsh1}-\eqref{pnbsh4+} and the upcoming Lemma \ref{pnbl1} we see
 \begin{eqnarray}\label{pnbresult}\mbfE\left(\log \varphi_n\right)\ra\mbfE\left(\log \varphi(D)\right).\end{eqnarray}
On the other hand, Lemma \ref{mto} and \eqref{au4} ensure $\mbfE\left(\log q_1\right)>-\infty$.
Another important fact is that \cite[Theorem 2.1]{lv1} told that
\begin{eqnarray}\lim_{D\ra\infty}\frac{\mbfE(\log\varphi(D))}{D}=-\frac{\gamma_{\sigma}}{\vartheta^2d^2}.\no\end{eqnarray}
Combining with \eqref{pnb0}, \eqref{pnbsh3} and recalling the definitions of $s_n$ and $\delta_n$ in \eqref{last62}, we see
\begin{eqnarray}\label{pnbsh6}&&\varliminf_{n\ra\infty}\mbfE\left(\frac{1}{n^{1/3}}\log P_n\right)
\geq -\frac{\gamma_{\sigma}}{\vartheta^2d^2}+2\delta\mbfE\left(\log q_1\right)\end{eqnarray} holds for $D$ large enough. According to \eqref{pnb3}, \eqref{pnb4} and \eqref{pnbsh6}, we see
$$\varliminf_{n\ra\infty} \mbfE\left(\frac{1}{n^{1/3}}\log p(n,b)\right)\geq-\frac{\gamma_{\sigma}}{\vartheta^2d^2}+2\delta\mbfE\left(\log q_1\right)-\vartheta (b-b'+d)+\frac{N-1}{N}\frac{\gamma_{\sigma}}{\vartheta^2d^2}.$$
We finally complete the proof by taking $N\ra+\infty, b'\ra b$ (and hence $\delta\ra 0$ because the choice of $\delta$), $d\ra \sqrt{\frac{\gamma_{\sigma}}{\vartheta^3 b}}$ (in this order).\qed
~~~
\begin{lemma}\label{pnbl1}
Let $g$ be a function with continuous and strictly positive derivative in its domain (a subset of $\bfR$) and the range of $g$ is $\bfR.$\footnote{We agree that $g(\cdot)=+\infty$ or $-\infty$ by a continuous extension if $\cdot$ is an edge of the domain of $g$. In addition, we note that this lemma is relative trivial if $g$ is bounded.} The range of random variables $X_n,$ $Y_n, X$ are all in the closure of the domain of $g.$ If

~~~~~{\rm (a)} $X_n\ra X$ in distribution,

~~~~~{\rm (b)} $\varlimsup_{n\ra+\infty}\bfE(|g(X_n)|^{z}+|g(Y_n)|^{z})<+\infty$ for some $z>1$ and $\bfE|g(X)|<+\infty,$

~~~~~{\rm (c)} $\exists a_n\ra 0$ such that $\bfP(|X_n-Y_n|\geq a_n)\ra 0$,

\noindent then $\bfE g(Y_n)\ra \bfE g(X).$
\end{lemma}
{\bf Proof of Lemma \ref{pnbl1}}: Since what we need to proof is about convergence, we can view the assumption (b) as $\sup_n\bfE(|g(X_n)|^{z}+|g(Y_n)|^{z})<+\infty,$ which implies that $\{g(X_n)\}$ and $\{g(Y_n)\}$ are both uniformly integrable, and hence for any $\ee>0,$ there exists $M$ such that
\begin{eqnarray}\label{uniint}\sup_n\left[\bfE(|g(X_n)|\1_{|g(X_n)|\geq M-1})+\bfE(|g(Y_n)|\1_{|g(Y_n)|\geq M-1})\right]+\bfE(|g(X)|\1_{|g(X)|\geq M})<\ee.\end{eqnarray}
Without loss of generality, we can also assume that $a_n$ satisfies $\bfP(|X_n-Y_n|\geq a_n)\leq a_n.$
Let $g^{-1}$ be the inverse function of $g$ and $l_n:=\sup_{x\in[g^{-1}(-M)-a_n, g^{-1}(M)+a_n]}|g'(x)|.$ For $n$ large enough, we have
\begin{eqnarray}\bfE(g(Y_n)\1_{|g(Y_n)|\leq M})
&=& \bfE(g(Y_n)\1_{|g(Y_n)|\leq M,|X_n-Y_n|\leq a_n})+\bfE(g(Y_n)\1_{|g(Y_n)|\leq M,|X_n-Y_n|> a_n})\no
\\&\leq&\bfE((g(X_n)+l_na_n)\1_{|g(Y_n)|\leq M,|X_n-Y_n|\leq a_n})+M\bfP(|X_n-Y_n|> a_n)\no
\\&\leq&\bfE(g(X_n)\1_{|g(Y_n)|\leq M,|X_n-Y_n|\leq a_n})+(l+M)a_n,\no\end{eqnarray}
where $l:=\varlimsup_{n\ra +\infty}l_n+1.$ We see $l<+\infty$ due to the definition of $l_n$. 
Furthermore, for $n$ large enough, it is true that
\begin{eqnarray*}
&&\bfE(g(X_n)\1_{|g(Y_n)|\leq M,|X_n-Y_n|\leq a_n})
\\&=&\bfE(g(X_n)\1_{g(X_n)>0,|g(Y_n)|\leq M,|X_n-Y_n|\leq a_n})+\bfE(g(X_n)\1_{g(X_n)<0,|g(Y_n)|\leq M,|X_n-Y_n|\leq a_n})
\\&\leq&\bfE(g(X_n)\1_{0<g(X_n)\leq M+la_n})+\bfE(g(X_n)\1_{la_n-M\leq g(X_n)<0,|X_n-Y_n|\leq a_n}),
\end{eqnarray*}
and
\begin{eqnarray*}
&&\bfE(g(X_n)\1_{la_n-M\leq g(X_n)<0,|X_n-Y_n|\leq a_n})
\\&=&\bfE(g(X_n)\1_{la_n-M\leq g(X_n)<0})-\bfE(g(X_n)\1_{la_n-M\leq g(X_n)<0,~|X_n-Y_n|> a_n})
\\&\leq&\bfE(g(X_n)\1_{la_n-M\leq g(X_n)<0})-(la_n-M)\bfE(\1_{la_n-M\leq g(X_n)<0,~|X_n-Y_n|> a_n})
\\&\leq&\bfE(g(X_n)\1_{la_n-M\leq g(X_n)<0})+(M-la_n)\bfP(|X_n-Y_n|> a_n).
\end{eqnarray*}
Therefore, it is true that
$$\bfE(g(Y_n)\1_{|g(Y_n)|\leq M})\leq\bfE(g(X_n)\1_{la_n-M\leq g(X_n)<M+la_n})+(2M+l)a_n.$$
By a similar argument, we can also obtain that
$$\bfE(g(Y_n)\1_{|g(Y_n)|\leq M})\geq\bfE(g(X_n)\1_{-la_n-M\leq g(X_n)<M-la_n})-(2M+l)a_n.$$
Hence for $n$ large enough, we have
\begin{eqnarray}&&|\bfE(g(Y_n)\1_{|g(Y_n)|\leq M})-\bfE(g(X_n)\1_{|g(X_n)|\leq M})|\no
\\&\leq&(2M+l)a_n+\bfE(|g(X_n)|\1_{|g(X_n)-M|\leq la_n})+\bfE(|g(X_n)|\1_{|g(X_n)+M|\leq la_n})\no
\\&\leq&(2M+l)a_n+\bfE(|g(X_n)|\1_{|g(X_n)|\geq M-1})\no\end{eqnarray}
and hence
\begin{eqnarray}\label{uniint2}
&&|\bfE g(Y_n)-\bfE g(X_n)|\no
\\&\leq&|\bfE(g(Y_n)\1_{|g(Y_n)|\leq M})-\bfE(g(X_n)\1_{|g(X_n)|\leq M})|+ \bfE(|g(X_n)|\1_{|g(X_n)|\geq M}+|g(Y_n)|\1_{|g(Y_n)|\geq M})\no
\\&\leq&(2M+l)a_n+2\ee.
\end{eqnarray}
From assumption (a) we have $|\bfE(g(X)\1_{|g(X)|\leq M})-\bfE(g(X_n)\1_{|g(X_n)|\leq M})|<\ee$ for $n$ large enough, hence \eqref{uniint} leads to the fact that $|\bfE g(X_n)-\bfE g(X)|<2\ee$ for $n$ large enough. Finally, combining the fact with \eqref{uniint2} we complete this proof. \qed

\begin{proposition}\label{pnbp2}
Under assumptions \eqref{au0} and {$\bf (A1)$}-{$\bf (A3)$},
it is true that
\begin{eqnarray}\label{pnb-ub}\varlimsup_{n\ra\infty}\mbfE(|n^{-1/3}\log p(n,b)|^{\nu})\leq +\infty\end{eqnarray} for some $\nu>\nu_0$.
\end{proposition}
{\bf Proof of Proposition \ref{pnbp2}} From \eqref{pnb2} one sees that it is enough to show
 $$\exists \nu>\nu_0,~~\varlimsup_{n\ra\infty}\mbfE(|n^{-1/3}\log P_n|^{\nu})< +\infty.$$
(Note that $\nu>1$.) Combining the inequality $|(\sum_{i=1}^nC_i)^{\nu}|\leq n^{\nu-1}\sum_{i=1}^n|C_i|^{\nu}$ with
\eqref{pnb>},
we see
 $$\mbfE(|\log P_n|^{\nu})\leq \delta^{\nu}_n\mbfE(|\log q_1|^{\nu})+2s_n^{\nu}\varlimsup_{n\ra\infty}[\mbfE(|\log \varphi_{n,end}|^{\nu})+\mbfE(|\log \varphi_n|^{\nu})]<+\infty.$$
Moreover, assumption \eqref{au4} and Lemma \ref{mto} mean that $\mbfE(|\log q_1|^{\lambda_6})<+\infty.$ Noting that $\delta_n=O(n^{1/3}), s_n=O(n^{1/3})$ and $\lambda_6>\nu_0$. Combining the above facts with \eqref{pnbsh3}, we see
\eqref{pnb-ub} holds for some $\nu\in(\nu_0, \lambda_6).$

\begin{proposition}\label{pnbp3} For any given $x>0$, there exists a constant $c_x>0$ such that
$$\mbfP\left(\frac{1}{n^{1/3}}\log p\left(n,b\right)<\frac{\gamma}{\sqrt{b}}-x\right)\leq c_xn^{(1-\nu)/3},$$
where $\nu$ is the one in Proposition \ref{pnbp2}.
\end{proposition}
Proof: We only need to consider the asymptotic behavior of $p(n,b)$ in the case $n=Nk$ due to the monotonicity of $p(n,b)$ (about $n$).
Recall the definition of $P_{lk,Nk}$ in \eqref{pnbpjn} and define
$$\bar{p}_{k,l,n}:=\sup_{y\in\bfR}\mbfP_{\L}\left(\forall_{(l-1)k\leq i\leq lk},\frac{T_i}{n^{1/3}}\in\left[\frac{b\vartheta i}{n}-d\vartheta ,\frac{b\vartheta i}{n}\right]\Big|T_{(l-1)k}=y\right),~~1\leq l\leq N,$$
then we see $P_{lk,Nk}\leq\Pi_{j=l+1}^{N}\bar{p}_{k,j,n}$ by Markov property.
On the other hand, define
$$\underline{p}_{k,1,n}:=\mbfP_\L
\left(I^*(0,k), T_{k}-c_{k,1,n}\in H_n\Big|T_{0}=0\right),$$
$$\underline{p}_{k,l,n}:=\inf_{z-c_{k,l-1,n}\in H_n}\mbfP_\L
\Big(I^*(lk-k,lk), T_{lk}-c_{k,l,n}\in H_n\Big|T_{(l-1)k}=z\Big), ~~2\leq l\leq N-1,$$
$$\underline{p}_{k,N,n}:=\inf_{z-c_{k,N-1,n}\in H_n}\mbfP_\L
\Big(I^*(Nk-k,Nk), T_{Nk}\geq \vartheta b'n^{\frac{1}{3}}\Big|T_{(N-1)k}=z\Big),$$
where $c_{k,l,n}:=\frac{\vartheta b(lk-\delta_n)}{n^{2/3}}, H_n:=[\lambda_5\delta_n,\delta_n/\lambda_5],$ ($\delta_n$ is the one defined in the proof of Proposition \ref{pnb1p})
then we have $P_{n}\geq\Pi_{j=1}^{N}\underline{p}_{k,j,n}$ by Markov property. Therefore, from \eqref{pnb3} we see
\begin{eqnarray*}p\left(n,b\right)\geq \frac{e^{n^{1/3}\vartheta b'}\Pi_{j=1}^{N}\underline{p}_{k,j,n}}{A_nk\sum\limits_{l=0}^{N-1}\left[e^{\vartheta\left(bn^{1/3}+dn^{1/3}-\frac{blk}{n^{2/3}}\right)}\Pi_{j=l+2}^{N}\bar{p}_{k,j,n}\right]},~~{\rm \mathbf{P}-a.s.}\end{eqnarray*}
where we agree that $\Pi_{j=N+1}^{N}\bar{p}_{k,j,n}=1.$ For any given $x$, we choose $d$ such that $\vartheta b>\frac{\gamma_\sigma}{\vartheta^2d^2}$ and $\vartheta d-\sqrt{\frac{\gamma_{\sigma}}{\vartheta b}}\leq\frac{x}{3}$. 
Denote
$$J_{k,l,n}:=\left\{\frac{\log \underline{p}_{k,l,n}}{n^{1/3}}\geq -\frac{\gamma_\sigma}{N\vartheta^2d^2}-\frac{x}{9N},~ \frac{\log \bar{p}_{k,l,n}}{n^{1/3}}\leq -\frac{\gamma_\sigma}{N\vartheta^2d^2}+\frac{x}{9N}\right\}.$$
In this proof, we take $A_n:=e^{\frac{x}{9}n^{1/3}}, b':=b-\frac{x}{9\vartheta}$. Then on $\cap_{l=1}^{N}J_{k,l,n},$ we have
\begin{eqnarray}\label{pnbp31}\frac{1}{n^{1/3}}\log p\left(n,b\right)&\geq& \vartheta (b'-b-d)-\frac{\gamma_\sigma}{\vartheta^2d^2}-\frac{x}{9}-\frac{x}{9}-\frac{1}{n^{1/3}}\log\left(k\sum_{l=0}^{N-1}e^{-\frac{\vartheta bl n^{1/3}}{N}-\frac{N-l-1}{N}(\frac{\gamma_\sigma}{\vartheta^2d^2}-\frac{x}{9})n^{1/3}}\right)\no
\\&\geq& -\vartheta d -\frac{\gamma_\sigma}{\vartheta^2d^2}-\frac{x}{3}-\frac{1}{n^{1/3}}\log\left(\frac{ke^{-\frac{N-1}{N}(\frac{\gamma_\sigma}{\vartheta^2d^2}-\frac{x}{9})n^{1/3}}}{1-e^{-\frac{\vartheta b n^{1/3}}{N}+\frac{n^{1/3}}{N}(\frac{\gamma_\sigma}{\vartheta^2d^2}-\frac{x}{9})}}\right).
\end{eqnarray}
Now we fix $N$ satisfying that $\frac{1}{N}\frac{\gamma_\sigma}{\vartheta^2d^2}<\frac{x}{9}$. Then for $k$ (recall $n=Nk$) large enough such that $\frac{\log (2k)}{n^{1/3}}<\frac{x}{9}$  and $\exp\{-\frac{\vartheta b n^{1/3}}{N}+\frac{n^{1/3}}{N}(\frac{\gamma_\sigma}{\vartheta^2d^2}-\frac{x}{9})\}<\frac{1}{2},$ \eqref{pnbp31} tells that
\begin{eqnarray}\frac{1}{n^{1/3}}\log p\left(n,b\right)&\geq&-\vartheta d -\frac{\gamma_\sigma}{\vartheta^2d^2}-\frac{x}{3}+\frac{N-1}{N}(\frac{\gamma_\sigma}{\vartheta^2d^2}-\frac{x}{9})-\frac{\log 2}{n^{1/3}},\no
\\&\geq& -\vartheta d-\frac{2x}{3}\no
\geq -\sqrt{\frac{\gamma_{\sigma}}{\vartheta b}}-x,~~\text{on}~\cap_{l=1}^{N}J_{k,l,n}.
\end{eqnarray}
Therefore, $\mbfP\left(\frac{1}{n^{1/3}}\log p\left(n,b\right)\geq-\sqrt{\frac{\gamma_{\sigma}}{\vartheta b}}-x\right)\geq \mbfP\left(\cap_{l=1}^{N}J_{k,l,n}\right).$
By the definition of $J_{k,l,n}$ one sees that $\{J_{k,l,n},l=1,2,...,N\}$ is an independent sequence of events and $\{\underline{p}_{k,l,n}, l=2,3,...N-1\}$ and $\{\bar{p}_{k,l,n},l=1,2,...N\}$ are both i.i.d. random sequences. Therefore,
\begin{eqnarray}\label{pnbp32}\mbfP\left(\frac{1}{n^{1/3}}\log p(n,b)\geq-\sqrt{\frac{\gamma_{\sigma}}{\vartheta b}}-x\right)\geq \mbfP\left(J_{k,1,n}\right)\mbfP\left(J_{k,2,n}\right)^{N-2}\mbfP\left(J_{k,N,n}\right).\end{eqnarray}
According to the forthcoming Lemma \ref{pnbl2}, we see that there exists a constant $c_5$ (depending on $x$, independent of $n$ and $k$) such that
\begin{eqnarray}\label{pnbp33}\mbfP\left(J^c_{k,2,n}\right)\leq c_5n^{\frac{1-\nu}{3}}~~\text{and}~~\mbfP\left(J^c_{k,N,n}\right)\leq c_5n^{\frac{1-\nu}{3}}.\end{eqnarray}
Recall that $q_m:=\mbfP_\L
(T_{m}\in [\lambda_5,1/\lambda_5],\xi_{m}\leq A_n|T_{m-1}=0)$ and $\delta_n:=\lfloor\delta n^{1/3}\rfloor$. For $J_{k,1,n},$ Markov property tells that 
$$\underline{p}_{k,1,n}\geq\prod_{i=1}^{\delta_n}q_i\inf_{z\in H_n}\mbfP_\L
\left(I(\delta_n,k), T_{k}-c_{k,1,n}\in H_n\Big|T_{\delta_n}=z\right),$$
According to Corollary \ref{strapp} and the facts $\mbfE(\log q_1)<0$ and $\mbfE(|\log q_1|^{\lambda_6})<+\infty,$ we see for $\delta$ small enough (for example, $0<\delta<\frac{-x}{2024N\mbfE(\log q_1)}$), there exists constant $c_6$ (independent of $x,n,k$) such that
$$\mbfP\left(\sum_{i=1}^{\delta_n}\log q_i\leq-\frac{x}{18N}n^{1/3}\right)\leq c_6n^{\frac{1-\lambda_6}{3}}\leq c_6n^{\frac{1-\nu}{3}}.$$
Moreover, Lemma \ref{pnbl2} also means that
$$\mbfP\left(\frac{\log\inf\limits_{z\in H_n}\mbfP_\L
\left(I(\delta_n,k), T_{k}-c_{n,k,1}\in H_n\Big|T_{\delta_n}=z\right)}{n^{1/3}}\leq-\frac{\gamma_\sigma}{N\vartheta^2d^2}-\frac{x}{18N}\right)\leq c_7n^{\frac{1-\nu}{3}}$$
holds for some $c_7$ (depending on $x$, independent of $n$ and $k$).
Therefore, combining with \eqref{pnbp33} we see there exists $c_8:=c_5+c_6+c_7$ such that
$\mbfP\left(J^c_{k,l,n}\right)\leq c_8n^{\frac{1-\nu}{3}}$ for any $l=1,2,...,N$. Recalling \eqref{pnbp32} we finally get
\begin{eqnarray}\label{pnbpzuihou}\mbfP\left(\frac{1}{n^{1/3}}\log\mbfP^{}_{\L}\left(n,\frac{b}{n^{2/3}}\right)\leq-\sqrt{\frac{\gamma_{\sigma}}{\vartheta b}}-x\right)\leq1-(1-c_8n^{\frac{1-\nu}{3}})^N.\end{eqnarray}
We recall that the choices of $c_8$ and $N$ are totally determined by $x$ (reviewing this proof, we only require $N$ to satisfy $\frac{1}{N}\frac{\gamma_\sigma}{\vartheta^2d^2}<\frac{x}{9}$), hence \eqref{pnbpzuihou} completes the proof.


\begin{lemma}\label{pnbl2}
For any given $b>0, d>0, r>\frac{\lambda_6}{\lambda_3(\lambda_6-\nu_0)}$, denote $$\bar{p}_{n}:=\sup_{y\in\bfR}\mbfP_{\L}\left(\forall_{i\leq n},\frac{T_i}{n^{1/3}}\in\left[\frac{b\vartheta i}{n}-d\vartheta ,\frac{b\vartheta i}{n}\right]\Big|T_{0}=y\right),$$
$$\underline{p}_{n}:=\inf_{z\in\left[x'n^{1/3},y'n^{1/3}\right]}\mbfP_\L
\left(\substack{\forall i\leq n,~~ T_i\in\left[\frac{\vartheta bi}{n^{2/3}}-\vartheta dn^{\frac{1}{3}},\frac{\vartheta bi}{n^{2/3}}\right],\\ T_{n}\in \left[x''n^{1/3},y''n^{1/3}\right],~~\xi_i\leq \exp\{n^{\frac{2}{3}r}\}}\Bigg|T_{0}=z\right),$$
where $-\vartheta d<x'< y'<0, \vartheta b-\vartheta d\leq x''< y''\leq \vartheta b.$
For any $\re>0$, there exists a constant $c_{\re}$ (depending only on $\re$) such that
\begin{eqnarray}\label{low}
\mbfP\left(\frac{\log\underline{p}_{n}}{n^{1/3}}<-\frac{\gamma_{\sigma}}{\vartheta^2d^2}-\re\right)\leq c_{\re}n^{(1-\nu)/3},
\end{eqnarray}
\begin{eqnarray}\label{up}
\mbfP\left(\frac{\log\bar{p}_{n}}{n^{1/3}}>-\frac{\gamma_{\sigma}}{\vartheta^2d^2}+\re\right)\leq c_{\re}n^{(1-\nu)/3},
\end{eqnarray}
where $\nu$ is the one in Proposition \ref{pnbp2}.
\end{lemma}
{\bf Proof of Proposition \ref{pnbl2}} Since the proof of \eqref{up} is similar and easier than the proof of \eqref{low}, here we only prove \eqref{low}.
In this proof we redefine some notations in the proof of Proposition \ref{pnb1p}. \footnote{Since in this lemma we set $-\vartheta d<x'$ and $y'<0$, a role like $\delta_n$ in Proposition \ref{pnb1p} is not necessary to be introduced.} Let
$$D_n=\lf Dn^{2/3}\rf,~~ D_{n,l}=l\lf Dn^{2/3}\rf,~~\text{and}~~ s_n:=\left\lf\frac{n}{\lf Dn^{2/3}\rf}\right\rf-1$$ and recall the definitions of $I(t_1, t_2)$ and $\left[c_1, c_2\right]_{n,l}$ in \eqref{pnbnotat}. Similar to the relationship in \eqref{pnb>}, for $n$ large enough, we have
\begin{eqnarray}
\underline{p}_{n}
\geq (\Pi_{l=1}^{s_n}\hat\varphi_{n,l})\hat\varphi_{n,end},
\end{eqnarray}
where
$$\hat\varphi_{n,l}:=\inf_{z\in\left[x'n^{1/3},y'n^{1/3}\right]}\mbfP_\L
\Big(I(D_{n,l-1},D_{n,l}), T_{D_{n,l}}\in \left[(\frac{x'+y'}{2})n^{1/3},y'n^{1/3}\right]\Bigg|T_{D_{n,l-1}}=z\Big),~$$
$$\hat\varphi_{n,end}:=\inf_{z\in\left[x'n^{1/3},y'n^{1/3}\right]}\mbfP_\L
\Big(I(D_{n,s_n},n), T_{n}\in \left[x'''n^{1/3},y'''n^{1/3}\right]\Bigg|T_{D_{n,s_n}}=z\Big)$$
and $x''':=\frac{2x''+y''}{3}, y''':=\frac{x''+2y''}{3}.$
Note that $\hat\varphi_{n,1},\hat\varphi_{n,2},...,\hat\varphi_{n,s_n},\hat\varphi_{n,end}$ are independent of each other and $\hat\varphi_{n,1},\hat\varphi_{n,2},...,\hat\varphi_{n,s_n}$ have the same distribution. Though an analog arguement in the proof of Proposition \ref{pnb1p} (the arguement on $\mbfE(\log\varphi_{n})$ and $\mbfE(\log\varphi_{n,end})$ therein),  we see that Lemma \ref{T-ubl} and Lemma \ref{T-ipl} imply that
\begin{eqnarray}\label{hatvarphi}&&\varlimsup\limits_{n\ra\infty}\mbfE(|\log\hat\varphi_{n,1}|^{\nu})+\varlimsup\limits_{n\ra\infty}\mbfE(|\log\hat\varphi_{n,end}|^{\nu})<+\infty~{\rm and}~\lim\limits_{n\ra\infty}\mbfE(\log\hat\varphi_{n,1})=\mbfE(\log\hat\varphi(D))~~~~~~~~\end{eqnarray} respectively, where
$$\hat\varphi(D):=\inf_{z\in\left[\frac{x'}{\sigma_*},\frac{y'}{\sigma_*}\right]}\mbfP_\L
\left(\begin{split}\forall s\leq D,~~ \sigma_*B_s+\sigma W_s\in\left[-\vartheta d,0\right],\\ \sigma_*B_D+\sigma W_D\in \left[\frac{x'+y'}{2},y'\right]\end{split}\Bigg|W,W_0=0,B_0=z\right)$$
and $B, W, \sigma, \sigma_*$  are the same ones in the proof of Proposition \ref{pnb1p}.
Moreover, \cite[Theorem 2.1]{lv1} tells that $\frac{\mbfE(\log\hat\varphi(D))}{D}\ra -\frac{\gamma_\sigma}{\vartheta^2d^2}.~$ 
Therefore, for any given $\re$, we can
choose $D$ large enough such that $\frac{\mbfE(\log\hat\varphi(D))}{D}\geq -\frac{\gamma_\sigma}{\vartheta^2d^2}-\frac{\re}{4}.$ Furthermore, for the choosen $D$, we have
$$\frac{s_n\mbfE(\log\hat\varphi_{n,1})+\mbfE(\log\hat\varphi_{n,end})}{n^{1/3}}\geq\frac{\mbfE(\log\hat\varphi(D))}{D}-\frac{\re}{4}$$
for $n$ large enough. From the above analysis, we see
\begin{eqnarray}&&\mbfP\left(\frac{\log\underline{p}_{n}}{n^{1/3}}<-\frac{\gamma_{\sigma}}{\vartheta^2d^2}-\re\right)\no
\\&\leq&\mbfP\left(\sum_{l=1}^{s_n}(\log\hat\varphi_{n,l}-\mbfE(\log\hat\varphi_{n,l}))+(\log\hat\varphi_{n,end}-\mbfE(\log\hat\varphi_{n,end}))<-\frac{\re}{2}n^{1/3}\right).\end{eqnarray}
We finally complete the proof by Corollary \ref{strapp} and \eqref{hatvarphi}.
\section{Preparation 4: A lower deviation for BPre}

In this section, let us forget the story of BRWre temporarily and focus on an asymptotic behavior of branching process in the random environment $\L$, \footnote{Since the motivation of this section is to consider some asymptotic behaviors of $p_2(l,\ee)$ in \eqref{p2} and the environment involved in \eqref{p2} is totally determined by $\L$, here we also denote $\L$ the random environment.} which is another important support to the proof in Section 8.
\begin{lemma}\label{7.1} Let $\{Z_n\}$ be a supercritical BPre (i.e., $\mbfE(\log \mbfE_{\L}Z_1)>0$)and
\begin{eqnarray}\label{BPreau2}\mbfE(|\min(0,\log \mbfE_{\L}Z_1)|^{\alpha})+\mbfE(|\max(0,\log \mbfE_{\L}(Z^2_1))|^{\alpha})<+\infty.\end{eqnarray}
for some $\alpha\geq 2.$ $g(x)$ is a positive-valued function satisfying $\lim_{x\ra\infty}g(x)=+\infty$. Let $b$ be a constant in $\left(1,e^{\mbfE(\log \mbfE_{\L}Z_1)}\right)$.
Denote $\eta_n:=\mbfP_{\L}(Z_n\geq \max\{b^{-n}\mbfE_{\L} Z_n,1\}).$
The following conclusions are true.

{\rm (1)} For any $l<\alpha$, $\varlimsup_{n\ra +\infty}\mbfE\left(|\frac{1}{g(n)}\log \eta_n|^l\right)<\infty$ when $\varlimsup_{n\ra +\infty}n^{2+l-\alpha}g(n)^{-l}<+\infty.$

{\rm (2)} For any given $x>0,$ $\varlimsup_{n\ra +\infty}g(n)^{\alpha-2}\mbfP\left(|\frac{1}{g(n)}\log \eta_n|>x\right)<+\infty.$

{\rm (3)} Suppose that $\alpha>2.$ We have $\lim_{n\ra +\infty}\mbfE\left(\frac{1}{g(n)}\log \eta_n\right)=0$ when $\varlimsup_{n\ra +\infty}n^{2+l-\alpha}g(n)^{-l}<+\infty$ holds for some $l>1$.

{\rm (4)} $\varlimsup_{n\ra +\infty}\mbfE\left(|\frac{1}{n}\log\max\{b^{-n}\mbfE_{\L} Z_n,1\}|^\alpha\right)\leq +\infty$.

{\rm (5)} $\varlimsup_{n\ra +\infty}n^{\alpha-1}\mbfP\left(\frac{1}{n}\log\max\{b^{-n}\mbfE_{\L} Z_n,1\}<\mbfE(\log \mbfE_{\L}Z_1)-\log b-x\right) <\infty$ for any given $x>0.$

{\rm (6)} $\varliminf_{n\ra +\infty}\mbfE\left(\frac{1}{n}\log\max\{b^{-n}\mbfE_{\L} Z_n,1\}\right)\geq \mbfE(\log \mbfE_{\L}Z_1)-\log b.$
\end{lemma}

{\bf Proof of Lemma \ref{7.1} (1)}
From the Paley-Zygmund inequality, for any $b>1$, we see
\begin{eqnarray}\label{PZine}\mbfP_{\L}(Z_{n}\geq b^{-n}\mbfE_{\L} Z_{n})\geq(1-b^{-n})^2\frac{(\mbfE_{\L} Z_{n})^2}{\mbfE_{\L} Z^2_{n}}, ~~{\rm \mathbf{P}-a.s.}\no\end{eqnarray}
Moreover, (note that $Z_n\in\bfN$,) Cauchy-Schwartz inequality tells that
\begin{eqnarray}\label{PZine}\mbfP_{\L}(Z_n\geq 1)\geq \frac{(\mbfE_{\L} Z_n)^2}{\mbfE_{\L} Z^2_n}, ~~{\rm \mathbf{P}-a.s.}\no\end{eqnarray}
Therefore, 
$$\eta_n\geq (1-b^{-n})^2\frac{(\mbfE_{\L} Z_{n})^2}{\mbfE_{\L} (Z^2_{n})}, ~~{\rm \mathbf{P}-a.s.}$$
In the rest of the proof we default all the equalities and inequalities hold in the sense of ${\rm \mathbf{P}-a.s.}$ unless we specially declare. 
Introduce the generate function of BPre as $f_{m,n}(s):=\sum_{i=0}^{+\infty}p_{m,n,i}s^{i},$ where $$p_{m,n,i}:=\mbfP_{\L}(Z_n=i|Z_m=1).$$
By some basic calculations and the relationship $f_{0,n}(s)=f_{0,n-1}(f_{n-1,n}(s))$, we see
\begin{eqnarray*}\frac{\mbfE_{\L}[Z_n(Z_n-1)]}{(\mbfE_{\L}Z_{n})^2}&=&\frac{f_{0,n}''(1)}{f_{0,n}'(1)^2}=\frac{f''_{0,n-1}(1)f'_{n-1,n}(1)^2+f_{0,n-1}'(1)f''_{n-1,n}(1)}{f_{0,n-1}'(1)^2f'_{n-1,n}(1)^2}
\\&=&\frac{f''_{0,n-1}(1)}{f_{0,n-1}'(1)^2}+\frac{f''_{n-1,n}(1)}{f_{0,n-1}'(1)f'_{n-1,n}(1)^2}.
\end{eqnarray*}
Iterating the above steps and note that $f_{0,0}(s)=s$, we get
\begin{eqnarray*}
\frac{\mbfE_{\L}[Z_n(Z_n-1)]}{(\mbfE_{\L}Z_{n})^2}&=&\sum_{i=1}^{n}\frac{f''_{i-1,i}(1)}{f_{0,i-1}'(1)f'_{i-1,i}(1)^2}.
\end{eqnarray*}
Hence for $n$ large enough, it is true that
\begin{eqnarray}\label{BPre1}
\frac{1}{\mbfP_{\L}(Z_n\geq \max\{b^{-n}\mbfE_{\L} Z_n,1\})}
\leq\frac{2\mbfE_{\L}(Z^2_n)}{(\mbfE_{\L}Z_{n})^2}=\sum_{i=1}^{n}\frac{2f''_{i-1,i}(1)}{f_{0,i-1}'(1)f'_{i-1,i}(1)^2}+\frac{2}{f_{0,n}'(1)}.
\end{eqnarray}
Recall that $b<e^{\mbfE(\log \mbfE_{\L}Z_1)}$ and hence there exists $q\in(1,+\infty)$ large enough such that $b<e^{\mbfE(\log \min( \mbfE_{\L}Z_1,q))}$.
Denote $$M_0:=1,~M_{i}:=\prod_{k=1}^i\min(q,f_{k}'(1)),~ \beta^*_i:=\frac{\max(f''_{i}(1),q^{-1})}{[\min(q,f'_{i}(1))]^2},~ f_i:=f_{i-1,i}.$$
Note that $f_{0,i}'(1)=\prod_{k=1}^if_{k}'(1)$ and the RHS of \eqref{BPre1} will not less than $1.$ Hence to prove Theorem \ref{7.1} (1), it is enough to show
\begin{eqnarray}\label{BPregoal}
\varlimsup_{n\ra+\infty}\mathbf E\left(\left[\frac{1}{g(n)}\log\left(\sum_{i=1}^{n}\frac{2\beta^*_i}{M_{i-1}}+\frac{2}{M_{n}}\right)\right]^l\right)<+\infty,\end{eqnarray}
or equivalently, to show
\begin{eqnarray}\label{BPreequiv}
\exists m_0\in\bfN,~ \varlimsup_{n\ra+\infty}\mathbf \sum_{m=m_0}^{+\infty}\mathbf P\left(\left[\frac{1}{g(n)}\log\left(\sum_{i=1}^{n}\frac{2\beta_i^*}{M_{i-1}}+\frac{2}{M_{n}}\right)\right]^l\geq m\right)<+\infty.\end{eqnarray}
Let us prove \eqref{BPreequiv} now.  Denote $\bar{b}:=\sum_{i=0}^{+\infty}b^{-i}<+\infty$ and $\beta_i:=2\beta_i^*+2$, we have
\begin{eqnarray}\label{defQ}
&&\mathbf P\left(\sum_{i=1}^{n}\frac{2\beta_i^*}{M_{i-1}}+\frac{2}{M_{n}}\geq e^{g(n)m^{1/l}}\right)\no
\\&\leq&\sum_{i=1}^{n+1}\mathbf P\left(\frac{\bar{b}\beta_i}{M_{i-1}}\geq e^{g(n)m^{1/l}}b^{1-i}\right)\no
\\&\leq&\sum_{i=0}^{n}\mathbf P\left(\log \beta_{i+1}-\mbfE\log \beta_1-\log M_i+i\mbfE\log M_1\geq g(n)m^{1/l}+i\mbfE\log M_1-i\log b-\mbfE\log (\beta_1\bar{b})\right)\no
\\&:=&\sum_{i=0}^{n}Q_{i,m}.\end{eqnarray}
Now we consider the upper bound of $Q_{i,m}$ based on the divergence rate of $g.$

From \eqref{BPreau2} we see $\mbfE(|\log M_1|^\alpha)+\mbfE(|\log \beta_1|^\alpha)<+\infty.$ Note that $\{f_i,i\in\bfN\}$ is an independent sequence, hence Corollary \ref{strapp} tells that we can find positive constants $c_9, c_{10}, c_{11}$ (all are independent of $n$) such that for $n$ large enough,
\begin{eqnarray}\label{BP1}
&&\sum_{m>(n/g(n))^l}\sum_{i=0}^{n}Q_{i,m}~~~~~~~({\rm note~that~}\log b<\mbfE\log M_1~{\rm and}~\mbfE\log (\beta_1\bar{b})<+\infty)\no
\\&\leq&\sum_{m>(n/g(n))^l}\sum_{i=0}^{n}\mathbf P\left(\log \beta_{i+1}-\mbfE(\log \beta_1)-\log M_i+i\mbfE(\log M_1)\geq c_9g(n)m^{1/l}\right)\no
\\&=&\sum_{m>(n/g(n))^l}\sum_{i=0}^{n}\mathbf P\left(\log \beta_{1}-\mbfE(\log \beta_1)-\log M_{i+1}+\log M_1+i\mbfE(\log M_1)\geq c_9g(n)m^{1/l}\right)\no
\\&\leq&\sum_{m>(n/g(n))^l}(n+1)\mathbf P\left(\max_{0\leq i\leq n}\big|\log \beta_1-\mbfE(\log \beta_1)-\log M_{i+1}+\log M_{1}+i\mbfE(\log M_1)\big|\geq c_9g(n)m^{1/l}\right)\no
\\&\leq&\sum_{m>(n/g(n))^l}\frac{c_{10}n^2}{[g(n)m^{1/l}]^{\alpha}}\leq c_{11}n^{2+l-\alpha}g(n)^{-l}~~~({\rm recall ~~that}~~ l<\alpha). 
\end{eqnarray}
(Throughout this paper we agree that the meaning of $\sum_{k>x}$ is $\sum_{k>x,k\in\bfN}$ if $x$ is not an integer.)
Therefore, if $\varlimsup_{n\ra +\infty} n/g(n)<+\infty,$ then we choose $m_0:=(\varlimsup_{n\ra +\infty} n/g(n))^l+1$ and thus \eqref{BPreequiv} follows from \eqref{BP1}.

Next we consider the case that $\varlimsup_{n\ra +\infty} n/g(n)=+\infty.$ 
Analog to the discussion in \eqref{BP1}, one can find constants $c_{12}$ and $c_{13}$ such that
\begin{eqnarray}\label{BP1+}
\sum_{m\leq (n/g(n))^l}~\sum_{i\leq g(n)m^{1/l}}Q_{i,m}
&\leq&\sum_{m\leq (n/g(n))^l}\frac{c_{12}[g(n)m^{1/l}]^2}{[g(n)m^{1/l}]^{\alpha}}\no
\\&\leq& c_{12}g(n)^{2-\alpha}\left(\sum_{m<(n/g(n))^l}m^{(2-\alpha)/l}\right)\no
\\&\leq& c_{13}g(n)^{2-\alpha}\left[\left(\frac{n}{g(n)}\right)^{l+2-\alpha}\1_{l>\alpha-2}+\log \left(\frac{n}{g(n)}\right)\1_{l=\alpha-2}+1\right].~~~~~~~
\end{eqnarray}
Since we have assumed that $\varlimsup_{n\ra +\infty}n^{2+l-\alpha}g(n)^{-l}<+\infty,$ the case $\varlimsup_{n\ra +\infty} n/g(n)=+\infty$ means $\alpha>2.$
Note that $f'_1(1)=\mbfE_{\L}Z_1$ and hence $\mbfE\log(M_1/b)>0$. From Corollary \ref{strapp} we can also find positive constants $c_{14}$ and $c_{15}$ such that
\begin{eqnarray}\label{BP1++}
&&\sum_{m\leq (n/g(n))^l}~\sum_{g(n)m^{1/l}<i\leq n}Q_{i,m}\no
\\&\leq&\sum_{m\leq (n/g(n))^l}~\sum_{g(n)m^{1/l}<i\leq n}\mathbf P\big(\log \beta_1-\mbfE(\log \beta_1)-\log M_{i+1}+\log M_{1}+i\mbfE(\log M_1)\geq i\mbfE\log(M_1/b)\big)\no
\\&\leq&\sum_{m\leq (n/g(n))^l}\sum_{g(n)m^{1/l}<i\leq n}c_{14}i^{1-\alpha}\no
\\&\leq&\sum_{m\leq (n/g(n))^l}c_{15}[g(n)m^{1/l}]^{2-\alpha}. ~~
\end{eqnarray}
holds for $n$ large enough. Note that the last line in \eqref{BP1++} is the same as the second line in  \eqref{BP1+} up to a multiplicative constant.
Combining \eqref{BP1} with \eqref{BP1+} and \eqref{BP1++}, we complete the proof of (1).

\noindent {\bf Proof of Lemma \ref{7.1} (2)}
We observe that if we let $l=1$ in the definition of $Q_{i,m}$ in \eqref{defQ}, then we have
$$\mbfP\left(\left|\frac{1}{g(n)}\log \eta_n\right|>x\right)\leq\sum_{i=0}^nQ_{i,x}.$$
From \eqref{BP1} we see $\sum_{i=0}^nQ_{i,x}$ will be dominated by $\frac{n^2}{[g(n)x]^{\alpha}}$ times a constant when $x>n/g(n);$ from \eqref{BP1+} and \eqref{BP1++} we see $\sum_{i=0}^nQ_{i,x}$ will be dominated by $[g(n)x]^{2-\alpha}$ times a constant when $x\leq n/g(n).$ Therefore, $\mbfP\left(|\frac{1}{g(n)}\log \eta_n|>x\right)$ will be dominated by $(g(n))^{2-\alpha}$ for any fixed $x$, which is the right conclusion in Lemma \ref{7.1} (2).

\noindent {\bf Proof of Lemma \ref{7.1} (3)}
To proof (3), it is enough to show $$\forall \ee>0,~~\lim_{n\ra +\infty}\mbfE\left(\frac{1}{g(n)}\log \eta_n\1_{\frac{1}{g(n)}\log \eta_n+\ee<0}\right)=0$$ for any given $\ee>0,$ which follows from Lemma \ref{7.1} (1), (2) and a standard application of Holder's inequality.

\noindent {\bf Proof of Lemma \ref{7.1} (4)} First we see $\mbfE\left(|\log f'_1(1)|^\alpha\right)\big(=\mbfE\left(|\log \mbfE_{\L} Z_1|^\alpha\right)\big)<+\infty$ due to \eqref{BPreau2} and the Jensen's inequality.
Note that
$$\mbfE\left(\left|\frac{1}{n}\log\max\{b^{-n}\mbfE_{\L} Z_n,1\}\right|^\alpha\right)=\mbfE\left(\left|\frac{\log\mbfE_{\L} Z_n}{n}-\log b\right|^\alpha\1_{\mbfE_{\L} Z_n>b^{n}}\right).$$
Hence it is enough to show
$$\varlimsup_{n\ra +\infty}\mbfE\left(\left|\frac{\log\mbfE_{\L} Z_n}{n}\right|^\alpha\right)<+\infty.$$
Recalling \eqref{BPreau2} we see the above inequality is obvious because of $\mbfE\left(|\log f'_1(1)|^\alpha\right)<+\infty$ and
 $$|\log\mbfE_{\L} Z_n|^{\alpha}= \left|\sum_{i=1}^n\log f'_i(1)\right|^{\alpha}\leq n^{\alpha-1}\sum_{i=1}^n\left|\log f'_i(1)\right|^\alpha.$$
{\bf Proof of Lemma \ref{7.1} (5)} Note that
\begin{eqnarray}&&\mbfP\left(\frac{1}{n}\log\max\{b^{-n}\mbfE_{\L} Z_n,1\}<\mbfE(\log \mbfE_{\L} Z_1)-\log b-x\right)\no
\\&\leq&\mbfP\left(\mbfE_{\L} Z_n<b^{n}\right)+\mbfP\left(\frac{\log\mbfE_{\L} Z_n}{n}<\mbfE(\log \mbfE_{\L} Z_1)-x\right).\no\end{eqnarray}
 Recall that $\mbfE\left(|\log \mbfE_{\L} Z_1|^\alpha\right)<+\infty$ and $\mbfE(\log \mbfE_{\L} Z_1)>\log b$. Hence Corollary \ref{strapp} tells that there exists constant $c_{16}$ (independent of $n$) such that $$\mbfP\left(\frac{1}{n}\log\max\{b^{-n}\mbfE_{\L} Z_n,1\}<\mbfE(\log \mbfE_{\L} Z_1)-\log b-x\right)\leq c_{16}n^{1-\alpha}.$$
\noindent{\bf Proof of Lemma \ref{7.1} (6)} In the proof of (5), we have mentioned that $$\varlimsup _{n\ra +\infty}n^{\alpha-1}\mbfP(\mbfE_{\L} Z_n\leq b^n)<+\infty~ {\rm and~ hence} ~\sum_{n=1}^{+\infty}\mbfP(\mbfE_{\L} Z_n\leq b^n)<+\infty.$$
Applying the Borel-Cantelli lemma we see $\mbfP\left(\lim_{n\ra+\infty}\frac{\log\max\{b^{-n}\mbfE_{\L} Z_n,1\}}{\log (b^{-n}\mbfE_{\L} Z_n)}=1\right)=1$.
Then the conclusion can be obtained by the law of large number and the fatou's lemma directly.

\section{Proof of \eqref{tpm}: the second step}

At the end of the Section 4, we say (to complete the proof of \eqref{tpm}) the only rest task is to find the upper bound of $\mbfP\left(\sqrt{\varepsilon}\log\left[\sum_{i=1}^{+\infty}\frac{\lc\bar m(i,\varepsilon)\rc}{\Pi_{j=1}^{i}f'_{j,\varepsilon}(1)}\right]>-\gamma+x\right)$ for a given $x.$ Now let us do it.
We note that the notation in this section is in line with the notation in Section 4 unless stated otherwise.

{\it Step 2} Note that $f'_{j,\varepsilon}(1)=p_1(j,\varepsilon)p_2(j,\varepsilon)\lc\bar m(j,\varepsilon)\rc$ and $\bar m(j,\varepsilon)\geq 1$, hence for $\ee>0$ small enough, we see
\begin{eqnarray}\label{st2-0}
&&\mbfP\left(\sqrt{\varepsilon}\log\left[\sum_{i=1}^{+\infty}\frac{\lc\bar m(i,\varepsilon)\rc}{\Pi_{j=1}^{i}f'_{j,\varepsilon}(1)}\right]>-\gamma+5x\right)\no
\\&\leq&\mbfP\left(\sqrt{\varepsilon}\log\left[\sum_{i=1}^{+\infty}\frac{2\bar m(i,\varepsilon)}{\Pi_{j=1}^{i}f'_{j,\varepsilon}(1)}\right]>-\gamma+5x\right)\no
\\&\leq&\mbfP\left(\sqrt{\varepsilon}\log\left[\sum_{i=1}^{+\infty}\frac{\bar m(i,\varepsilon)}{\Pi_{j=1}^{i}f'_{j,\varepsilon}(1)}\right]>-\gamma+4x\right)\no
\\&\leq&\mbfP\left(\sqrt{\varepsilon}\log\frac{1}{p_1(1,\varepsilon)p_2(1,\varepsilon)}>-\gamma+2x\right)+\mbfP\left(\sqrt{\varepsilon}\log\left[1+\sum_{i=2}^{+\infty}\frac{\bar m(i,\varepsilon)}{\bar m(1,\varepsilon)\Pi_{j=2}^{i}f'_{j,\varepsilon}(1)}\right]>2x\right)\no
\\&:=&\psi_{\ee}+\psi^*_{\ee}.
\end{eqnarray}
Let $v$ be a positive constant and hence $\sum_{i=1}^{+\infty}e^{\frac{-iv}{\sqrt{\ee}}}\leq2<3\leq e^{\frac{x}{\sqrt{\ee}}}$ holds for $\ee$ small enough. Then we see
\begin{eqnarray}\label{st2-1}
\psi^*_{\ee}&\leq&
\mbfP(1>e^{\frac{x}{\sqrt{\ee}}})+\mbfP\left(\sum_{i=2}^{+\infty}\frac{\bar m(i,\varepsilon)}{\bar m(1,\varepsilon)\Pi_{j=2}^{i}f'_{j,\varepsilon}(1)}>2e^{\frac{x}{\sqrt{\ee}}}\right)~~~({\rm note~that}~e^{\frac{2x}{\sqrt{\ee}}}>3e^{\frac{x}{\sqrt{\ee}}})\no
\\&\leq&
\sum_{i=2}^{+\infty}\mbfP\left(\frac{\bar m(i,\varepsilon)}{\bar m(1,\varepsilon)\Pi_{j=2}^{i}f'_{j,\varepsilon}(1)}>e^{\frac{x-(i-1)v}{\sqrt{\ee}}}\right)\no
\\&=&\sum_{i=2}^{+\infty}\mbfP\left(\sqrt{\ee}\log \bar{m}(1,\varepsilon)+\sqrt{\ee}\log\Pi_{j=2}^{i-1}f'_{j,\varepsilon}(1)+\sqrt{\ee}\log(p_1(i,\varepsilon)p_2(i,\varepsilon))<(i-1)v-x\right)~~~\no
\\&\leq&\sum_{i=2}^{+\infty}\mbfP\left(\sqrt{\ee}\log \bar{m}(1,\varepsilon)+\sqrt{\ee}\log\Pi_{j=2}^{i-1}f^*_{j,\varepsilon}+\sqrt{\ee}\log(p_1(i,\varepsilon)p_2(i,\varepsilon))<(i-1)v\right)\no
\\&:=&\sum_{i=2}^{+\infty}\psi^*_{\ee,i},
\end{eqnarray} where $f^*_{j,\varepsilon}:=p_1(j,\varepsilon)p_2(j,\varepsilon)\bar m(j,\varepsilon)\leq f'_{j,\varepsilon}(1)$ and we agree that $\Pi_{j=n+1}^{n}\cdot=1.$
Denote $\mu_{\ee}:=\mbfE(\sqrt{\ee}\log f^*_{1,\varepsilon}).$ Recall the relationship $n=\left\lfloor(\varsigma+z)\varepsilon^{\frac{-3}{2}}\right\rfloor$ in \eqref{fhao} and thus $\lim_{\ee\downarrow 0}\sqrt{\ee}(\varsigma n)^{1/3}=\varsigma^{2/3}$, which means that for any $c'<1$, it is true that $\ee i\geq \frac{c'\varsigma^{4/3}i}{(\varsigma n)^{2/3}}$ 
as long as $\ee$ small enough. 
Applying Proposition \ref{pnb1p} we obtain that \begin{eqnarray}\label{st-p11}\varliminf_{\ee\downarrow 0}\sqrt{\ee}\mbfE(\log p_1(1,\varepsilon))\geq -\varsigma^{2/3}\sqrt{\frac{\gamma_{\sigma}}{c'\varsigma^{4/3}\vartheta}}=-\sqrt{\frac{\gamma_{\sigma}}{c'\vartheta}}~{\rm and~thus}~\varliminf_{\ee\downarrow 0}\sqrt{\ee}\mbfE(\log p_1(1,\varepsilon))\geq \gamma.~~\end{eqnarray}
Recall the construction of $\{Z_i(l)\}_{i\in\bfN}$ in \eqref{defzn} and the definition of $p_2(1,\varepsilon)$ in \eqref{p2}. From the assumptions \eqref{au3} and \eqref{au6} and the notation $\lambda_8:=\min(\lambda_3,\lambda_7)$ one sees that
\begin{eqnarray}\label{auzn}\forall \lambda\in(2,\lambda_8),~~~ l\in\bfN^+,~~~\mbfE(|\log^- \mbfE_{\L}Z_1(l)|^{\lambda})+\mbfE(|\log^+ \mbfE_{\L}(Z^2_1(l))|^{\lambda})<+\infty.\end{eqnarray}
Then combining with Lemma \ref{7.1} (3), we get
\begin{eqnarray}\label{st-p12}\lim_{\ee\downarrow 0}\sqrt{\ee}\mbfE(\log p_2(1,\varepsilon))=0.\end{eqnarray} Note that $\lim_{\ee\downarrow 0}\sqrt{\ee}\lf zn\rf=\frac{(c-1)\varsigma^{2}}{a}$ and hence Lemma \ref{7.1} (6) means that \begin{eqnarray}\label{st-p13}\lim_{\ee\downarrow 0}\sqrt{\ee}\mbfE(\log\bar m(1,\varepsilon))\geq\frac{(c-1)\varsigma^{2}}{a}(m(a)-\log w).\end{eqnarray} From the above calculations we see $$\varliminf_{\ee\downarrow 0}\mu_{\ee}\geq \frac{(c-1)\varsigma^{2}}{a}(m(a)-\log w)+\gamma:=\underline\mu.$$ Recall that $m(a)>0$ and we have chosen 
$w\in(1,e^{m(a)}).$ Then for any given $c>1,$ we can find a $\varsigma\in\bfN$ large enough such that $\underline\mu>0$ and from now on, we take $v:=\frac{1}{2}\underline\mu$.

We next estimate the upper bound of $\psi^*_{\ee,i}$ in two different strategies.

\emph{Strategy 1} By the subadditivity of probability and the relationship $$v:=\frac{1}{2}\underline{\mu}=\left(\frac{(c-1)\varsigma^{2}}{a}(m(a)-\log w)-\frac{1}{6}\underline{\mu}\right)+\left(0-\frac{1}{6}\underline{\mu}\right)+\left(\gamma-\frac{1}{6}\underline{\mu}\right),$$ we have
\begin{eqnarray}\label{me0}\psi^*_{\ee,i}&\leq&\mbfP\left(\sqrt{\ee}\log(m(1,\varepsilon)p_1(i,\varepsilon)p_2(i,\varepsilon)) <\frac{1}{2}\underline{\mu}\right)+\mbfP\left(\sum_{j=2}^{i-1}\sqrt{\ee}\log f^*_{j,\varepsilon}<\frac{i-2}{2}\underline{\mu}\right)\no
\\&\leq&(i-1)\mbfP\left(\sqrt{\ee}\log p_1(1,\varepsilon)<\gamma-\frac{1}{6}\underline\mu\right)\no
+(i-1)\mbfP\left(\sqrt{\ee}\log p_2(1,\varepsilon)<-\frac{1}{6}\underline\mu\right)\no
\\&~~~&~~+(i-1)\mbfP\left(\sqrt{\ee}\log \bar{m}(1,\varepsilon)<\frac{(c-1)\varsigma^{2}}{a}(m(a)-\log w)-\frac{1}{6}\underline\mu\right).
\end{eqnarray}
Recalling that $\underline\mu>0$, Proposition \ref{pnbp3} tells that
\begin{eqnarray}\label{-p131}\exists \nu>\nu_0,~~~\varlimsup_{\ee\downarrow0}\ee^{\frac{1-\nu}{2}}\mbfP\left(\sqrt{\ee}\log p_1(1,\varepsilon)<\gamma-\frac{1}{6}\underline\mu\right)<+\infty.\end{eqnarray}
From \eqref{auzn} and Lemma \ref{7.1} (2) and (5), we have 
\begin{eqnarray}\label{-p231}~\varlimsup_{\ee\downarrow0}\ee^{\frac{2-\lambda_8}{2}}\mbfP\left(\sqrt{\ee}\log p_2(1,\varepsilon)<-\frac{1}{6}\underline\mu\right)<+\infty\end{eqnarray}
and
\begin{eqnarray}\label{-p33}\varlimsup_{\ee\downarrow0}\ee^{\frac{1-\lambda_8}{2}}\mbfP\left(\sqrt{\ee}\log m(1,\varepsilon)<\frac{(c-1)\varsigma^{2}}{a}(m(a)-\log w)-\frac{1}{6}\underline\mu\right)<+\infty.\end{eqnarray}
Combining \eqref{me0} with \eqref{-p131}-\eqref{-p33}, we can find a constant $c_{17}$ (independent of $i$) such that 
\begin{eqnarray}\label{me1}\forall i\geq 2,~~~ \psi^*_{\ee,i}<c_{17}i\ee^{\frac{1}{2}\min(1+\nu,\lambda_8)-1}\end{eqnarray}
holds for $\ee>0$ small enough.

\emph{Strategy 2} The another strategy to estimate the upper bound of $\psi^*_{\ee,i}$ is using Corollary \ref{strapp}. We observe that
\begin{eqnarray}\psi^*_{\ee,i}=\mbfP\left(\sqrt{\ee}\log (p_1(i,\varepsilon)p_2(i,\varepsilon)\bar{m}(1,\varepsilon))-\mu_{\ee}+\sum_{j=2}^{i-1}(\sqrt{\ee}\log f^{*}_{j,\varepsilon}-\mu_{\ee}))<(i-1)(v-\mu_{\ee})\right).\no \end{eqnarray}
Recall that $0<\underline{\mu}\leq \varliminf_{\ee\downarrow0}\mu_{\ee}$ and $v:=\frac{1}{2}\underline{\mu}$. Then for $\ee$ small enough, we have
\begin{eqnarray}\psi^*_{\ee,i}\leq\mbfP\left(\sqrt{\ee}\log (p_1(i,\varepsilon)p_2(i,\varepsilon)\bar{m}(1,\varepsilon))-\mu_{\ee}+\sum_{j=2}^{i-1}(\sqrt{\ee}\log f^*_{j,\varepsilon}-\mu_{\ee})<(1-i)\frac{1}{3}\underline{\mu}\right).\no \end{eqnarray}
Note that $\sqrt{\ee}\log(\bar{m}(1,\varepsilon)p_1(i,\varepsilon)p_2(i,\varepsilon))$ is independent of $\sqrt{\ee}\log f^*_{j,\varepsilon}, j=2,3,...,i-1$ and they all have the same expectation $\mu_{\ee}.$
Now we focus on the integrability of $\sqrt{\ee}\log f^*_{j,\varepsilon}.$ Proposition \ref{pnbl2} tells that
\begin{eqnarray}\exists \nu>\nu_0, ~~\label{-p12}\varlimsup_{\ee\downarrow0}\mbfE\left(|\sqrt{\ee}\log p_1(1,\varepsilon)|^{\nu}\right)<+\infty.\end{eqnarray}
Combining with \eqref{auzn}, Lemma \ref{7.1} (1) and (4) tell that
$$\forall \lambda<\lambda_8,~~\varlimsup_{\ee\downarrow0}\big[\mbfE(|\sqrt{\ee}\log p_2(1,\varepsilon)|^{\lambda})+\mbfE(|\sqrt{\ee}\log \bar m(1,\varepsilon)|^{\lambda_8})\big]<+\infty.$$
Then Corollary \ref{strapp} tells that for $\ee$ small enough, there exists a constant $c_{18}$ (independent of $i$) such that
\begin{eqnarray}\label{me2}\forall i\geq 2,~~~\psi^*_{\ee,i}<c_{18}i^{1-\min(\nu,\lambda)}.\end{eqnarray}
Moreover, recalling the definition of $\psi_{\ee}$ in \eqref{st2-0}, we see
\begin{eqnarray}\label{st2-2}\varlimsup_{\ee\downarrow0}\ee^{1-\frac{1}{2}\min(1+\nu,\lambda_8)}\psi_{\ee}<+\infty\end{eqnarray}
because \eqref{-p131}, \eqref{-p231} still hold if $\frac{1}{6}\underline{\mu}$ therein is replaced by $x$.

Note that we can choose $\lambda$ such that $\min(\nu,\lambda)>2$ due to $\lambda_8>2$ and $\nu_0>2$. Then for any positive-valued function $h$, \eqref{st2-1}, \eqref{me1}, \eqref{me2} and \eqref{st2-2} ensure that
\begin{eqnarray}\label{me3}\exists c_{19}, c_{20},~~~ \psi_{\ee}+\psi^*_{\ee}&\leq&\sum_{i<h(\ee)}c_{19}i\ee^{\frac{1}{2}\min(1+\nu,\lambda_8)-1}+\sum_{i\geq h(\ee)}c_{18}i^{1-\min(\nu,\lambda)}\no
\\&\leq& c_{20}h(\ee)^2\ee^{\frac{1}{2}\min(1+\nu,\lambda_8)-1}+c_{20}h(\ee)^{2-\min(\nu,\lambda)}.
\end{eqnarray}
Then choosing $h(\ee)$ such that $h(\ee)^{-\min(\nu,\lambda)}=\ee^{\frac{1}{2}\min(1+\nu,\lambda_8)-1}$ and recalling \eqref{st2-0}, we finally get \eqref{tpm} and thus $\lim_{\ee\downarrow 0}\mbfP(\sqrt{\ee}\log\varrho_{\L}(\varepsilon)<\gamma-x)=0$ holds for any given $x>0.$\qed
\section{Proofs of Theorem \ref{main} (4)-(6)}

\noindent{\bf Proof of Theorem \ref{main} (4)} Recalling the definition of $\ee_n$ in \eqref{fhao}, we can see $$\lim\limits_{n\ra+\infty}\sqrt{\ee_{n}}\log\varrho_{\L}(\varepsilon_{n})\geq \gamma, \mathbf{P}-\text{a.s.}$$ follows from Borel-Cantelli lemma and \eqref{tpm}. Therefore, we finish the proof of (4) by the relationship $\ee^{-1/2}_{n}=O(n^{1/3})$ and the fact that $\sqrt{\ee}\log\varrho_{\L}(\varepsilon)\geq \sqrt{\ee_{n-1}}\log\varrho_{\L}(\varepsilon_{n})$ when $\ee_{n}<\ee\leq\ee_{n-1}$.

\noindent{\bf Proof of Theorem \ref{main} (5)}
Note that $\mbfE(|\sqrt{\ee}\log\varrho_{\L}(\varepsilon)|^p)\no
\leq\sum_{k=0}^{+\infty}\mbfP(|\sqrt{\ee}\log\varrho_{\L}(\varepsilon)|^p\geq k).$ Hence according to \eqref{lowrho} and the trick used in \eqref{st2-0}, it is enough to show
\begin{eqnarray}\label{5-0}\bar\psi_{\ee}:=\sum_{k\geq|7\gamma|^{p}}^{+\infty}\mbfP\left(\sqrt{\varepsilon}\log\left[\sum_{i=1}^{+\infty}\frac{\bar m(i,\varepsilon)}{\Pi_{j=1}^{i}f^*_{j,\varepsilon}}\right]\geq k^{1/p}\right)<+\infty.\end{eqnarray}
For $\ee$ small enough such that $e^{\frac{|7\gamma|}{2\sqrt{\ee}}}\geq(1-e^{-\frac{v}{\sqrt{\ee}}})^{-1}$, (where $v:=\frac{1}{2}\underline{\mu}$ is has been introduced in Section 8,) it is true that
\begin{eqnarray}\label{5-1}
\bar\psi_{\ee}&\leq&\sum_{k\geq|7\gamma|^{p}}^{+\infty}\mbfP\left(\sum_{i=1}^{+\infty}\frac{\bar m(i,\varepsilon)}{\Pi_{j=1}^{i}f^*_{j,\varepsilon}}\geq e^{\frac{k^{1/p}}{2\sqrt{\ee}}}\sum_{i=1}^{+\infty}e^{-\frac{(i-1)v}{\sqrt{\ee}}}\right)\no
\\&\leq&\sum_{i=1}^{+\infty}\sum_{k\geq|7\gamma|^{p}}^{+\infty}\mbfP\left(\frac{\bar m(i,\varepsilon)}{\Pi_{j=1}^{i}f^*_{j,\varepsilon}}\geq e^{\frac{k^{1/p}}{2\sqrt{\ee}}}e^{-\frac{(i-1)v}{\sqrt{\ee}}}\right)\no
\\&\leq&\sum_{i=0}^{+\infty}\sum_{k\geq|7\gamma|^{p}}^{+\infty}\left[\mbfP\left(\sqrt{\ee}\log p_1(i+1,\varepsilon)p_2(i+1,\varepsilon)+\sum_{j=1}^{i}\sqrt{\ee}\log f^*_{j,\varepsilon}\leq iv-\frac{k^{1/p}}{2}\right)\right]\no
\\&:=&\sum_{i=0}^{+\infty}\sum_{k\geq|7\gamma|^{p}}^{+\infty}\psi(\ee,i,k).
\end{eqnarray}
Note that for $\ee$ small enough, we have $$iv-\frac{k^{1/p}}{2}-i\mu_{\ee}-\mbfE(\sqrt{\ee}\log p_1(1,\varepsilon)p_2(1,\varepsilon))\leq -\frac{1}{3}\underline{\mu}i-\frac{k^{1/p}}{3}.$$
Hence for $\ee$ small enough, we can find $c_{21}$ (independent of $i,k,\ee$) such that
\begin{eqnarray}\label{5-2}\forall \lambda<\lambda_8, ~~~\psi(\ee,i,k)\leq c_{21}i(\underline{\mu}i+k^{1/p})^{-\min(\nu,\lambda)}\end{eqnarray} by using Corollary \ref{strapp}.
On the other hand, (recall that $\gamma:=-\sqrt{\frac{\gamma_{\sigma}}{\vartheta}}$ and) note that
$$\psi(\ee,i,k)\leq \mbfP\left(\sqrt{\ee}\log p_1(i+1,\varepsilon)p_2(i+1,\varepsilon)+\sum_{j=1}^{i}\sqrt{\ee}\log f^*_{j,\varepsilon}\leq iv+2\gamma\right).$$
By an analogue of strategy 2 in Section 8, we see there exists a constant $c_{22}$ (independent of $i$) such that
\begin{eqnarray}\label{5-3}\psi(\ee,i,k)\leq c_{22}i\ee^{\frac{\min(\nu+1,\lambda_8)}{2}-1}.\end{eqnarray}
The assumption $2+p<\min(\nu_0,\lambda_8)$ means that we can find $\lambda<\lambda_8$ to satisfy $2+p<\min(\nu,\lambda)$. From \eqref{5-2} and \eqref{5-3} we see for any real function $h^*:[0,+\infty)\ra[\frac{|4\gamma|}{7v}, +\infty)$, there exist constants $c_{23}$ and $c_{24}$ such that
\begin{eqnarray}&&\sum_{i\geq\frac{|\gamma|}{v}}\sum_{|7\gamma|^{p}\leq k\leq (7vi)^{p}}\psi(\ee,i,k)\no
\\&\leq&c_{23}\sum_{\frac{|\gamma|}{v}\leq i\leq h^*(\ee)}i^{p+1}\ee^{\frac{\min(\nu+1,\lambda_8)}{2}-1}+c_{23}\sum_{i>h^*(\ee)}\sum_{|7\gamma|^{p}\leq k\leq (7vi)^{p}}(\underline{\mu}i+k^{1/p})^{1-\min(\nu,\lambda)}\no
\\&\leq&c_{24}\sum_{\frac{|\gamma|}{v}\leq i\leq h^*(\ee)}i^{p+1}\ee^{\frac{\min(\nu+1,\lambda_8)}{2}-1}+c_{24}\sum_{i>h^*(\ee)}i^{p+1-\min(\nu,\lambda)}\no
\\&\leq&c_{25}h^*(\ee)^{2+p}\ee^{\frac{\min(\nu+1,\lambda_8)}{2}-1}+c_{25}h^*(\ee)^{2+p-\min(\nu,\lambda)}.\no\end{eqnarray}
Since $\frac{\min(\nu+1,\lambda_8)}{2}-1\geq\frac{\min(\nu,\lambda_8)}{2}-1>0$,
it is easy to find a proper $h^*$ such that \begin{eqnarray}\label{5-4}\sum_{i\geq \frac{|4\gamma|}{7v}}\sum_{~|4\gamma|^{p}\leq k\leq (7vi)^{p}}\psi(\ee,i,k)<+\infty.\end{eqnarray}
Moreover, \eqref{5-2} also implies that we can find a constant $c_{26}$ such that
\begin{eqnarray}\sum_{i=0}^{+\infty}\sum_{k>(7vi)^{p}}\psi(\ee,i,k)\no
&\leq& c_{21}\sum_{i=0}^{+\infty}\sum_{k>(7vi)^{p}}ik^{\frac{-\min(\nu,\lambda)}{p}}
\leq c_{26}\sum_{i=0}^{+\infty}i^{1+p-\min(\nu,\lambda)}\no\end{eqnarray}
due to $\min(\nu,\lambda)>p$.
Therefore, the assumption $2+p<\min(\nu_0,\lambda_8)$ means that
\begin{eqnarray}\label{5-5}\sum_{i=0}^{+\infty}\sum_{k>(7vi)^{p}}\psi(\ee,i,k)<+\infty.\end{eqnarray}
Finally, combining the following equalities
\begin{eqnarray}
\sum_{i=0}^{+\infty}\sum_{k>|7\gamma|^{p}}^{+\infty}\psi(\ee,i,k)\no
&\leq&\left(\sum_{i=0}^{+\infty}\sum_{k>(7vi)^{p}}+\sum_{i\geq\frac{|\gamma|}{v}}\sum_{|7\gamma|^{p}\leq k\leq (7vi)^{p}}
+\sum_{i<\frac{|\gamma|}{v}}\sum_{|7\gamma|^{p}\leq k\leq (7vi)^{p}}
\right)\psi(\ee,i,k)\no
\\&=&\left(\sum_{i=0}^{+\infty}\sum_{k> (7vi)^{p}}+\sum_{i\geq\frac{|\gamma|}{v}}\sum_{|7\gamma|^{p}\leq k\leq (7vi)^{p}}
\right)\psi(\ee,i,k)+0\no
\end{eqnarray}
with \eqref{5-1}, \eqref{5-4} and \eqref{5-5}, we complete the proof of (5).

\noindent{\bf Proof of Theorem \ref{main} (6)} Under the assumptions in (6) we see $\sqrt{\ee}\log\varrho_{\L}(\varepsilon)$ converges to $\gamma$ in probability. Combining with the conclusion in (5), we see $\lim\limits_{\ee\downarrow 0}\mathbf{E}\left(|\sqrt{\ee}\log\varrho_{\L}(\varepsilon)-\gamma|^p\right)=0$ by a standard application of Holder's inequality. \qed


\end{document}